\pdfoutput=1
\RequirePackage{ifpdf}
\ifpdf 
\documentclass[pdftex]{sigma}
\else
\documentclass{sigma}
\fi

\numberwithin{equation}{section}
\newtheorem{teo}{Theorem}[section]

\theoremstyle{definition}
\newtheorem{dfn}[teo]{Definition}
\newtheorem{exmpl}[teo]{Example}
\newtheorem{rmk}[teo]{Remark}
\newtheorem{assumption}[teo]{Assumption}

\DeclareMathOperator{\tr}{tr}
\begin{document}

\allowdisplaybreaks

\newcommand{\arXivNumber}{1408.4088}

\renewcommand{\PaperNumber}{001}

\FirstPageHeading

\ShortArticleName{Geometry of Centroaf\/f\/ine Surfaces in $\mathbb{R}^5$}

\ArticleName{Geometry of Centroaf\/f\/ine Surfaces in $\boldsymbol{\mathbb{R}^5}$}

\Author{Nathaniel BUSHEK~$^\dag$ and Jeanne N.~CLELLAND~$^\ddag$}

\AuthorNameForHeading{N.~Bushek and J.N.~Clelland}

\Address{$^\dag$~Department of Mathematics, UNC - Chapel Hill,\\
\hphantom{$^\dag$}~CB \#3250, Phillips Hall, Chapel Hill, NC 27599, USA}
\EmailD{\href{mailto:bushek@unc.edu}{bushek@unc.edu}}

\Address{$^\ddag$~Department of Mathematics, 395 UCB, University of Colorado, Boulder, CO 80309-0395, USA}
\EmailD{\href{mailto:Jeanne.Clelland@colorado.edu}{Jeanne.Clelland@colorado.edu}}

\ArticleDates{Received August 23, 2014, in f\/inal form December 26, 2014; Published online January 06, 2015}

\Abstract{We use Cartan's method of moving frames to compute a~complete set of local invariants for nondegenerate,
2-dimensional centroaf\/f\/ine surfaces in $\mathbb{R}^5 \setminus \{0\}$ with nondegenerate centroaf\/f\/ine metric.
We then give a~complete classif\/ication of all homogeneous centroaf\/f\/ine surfaces in this class.}


\Classification{53A15; 58A15}

\section{Introduction}

An immersion $\bar{f}:M \to \mathbb{R}^n \setminus \{0 \}$ is called a~{\em centroaffine immersion}
(cf.~Def\/inition~\ref{centroaffine-immersion-def}) if the position vector $\bar{f}(x)$ is transversal to the tangent
space $\bar{f}_*(T_xM)$ for all $x\in M$.
{\em Centroaffine geometry} is the study of those properties of centroaf\/f\/ine immersions that are invariant under the
action of the {\em centroaffine group} ${\rm GL}(n,\mathbb{R})$ on $\mathbb{R}^n \setminus \{0 \}$.
Much attention has been given to the study of centroaf\/f\/ine curves and hypersurfaces~\cite{GW97, Laugwitz65, Wang91,Wang95,MM33,
Scharlach97, Scharlach96, Wang94}, and more recently to the study of centroaf\/f\/ine immersions of
codimension~2~\cite{Furuhata00, Furuhata02,NS93, Scharlach95,Scharlach99, Scharlach98, Yang08, Yang13}.
In~\cite{Wilkens99}, applications of centroaf\/f\/ine geometry to the problem of feedback equivalence in control theory are
discussed.

In this paper, we consider the case of a~2-dimensional centroaf\/f\/ine surface in $\mathbb{R}^5 \setminus \{0 \}$.
We will use Cartan's method of moving frames to construct a~complete set of local invariants for a~large class of such
surfaces under certain nondegeneracy assumptions (cf.\ Def\/inition~\ref{nondegenerate-def},
Assumption~\ref{constant-type-assumption}).
In addition, we give a~complete classif\/ication of the {\em homogeneous} centroaf\/f\/ine surfaces in this class~-- i.e.,
those that admit a~3-dimensional Lie group of symmetries that acts transitively on an adapted frame bundle canonically
associated to the surface (cf.\ Def\/inition~\ref{homogeneous-def}).
Our primary results are Theorems~\ref{space-like-thm} and~\ref{time-like-thm}, which describe local invariants for
centroaf\/f\/ine surfaces with def\/inite and indef\/inite centroaf\/f\/ine metrics, respectively, and Theorem~\ref{homog-thm},
which describes the homogeneous examples.

The paper is organized as follows.
In Section~\ref{background-sec}, we introduce the basic concepts of centroaf\/f\/ine geometry and centroaf\/f\/ine surfaces
in~$\mathbb{R}^5 \setminus \{0 \}$, including the centroaf\/f\/ine frame bundle and the Maurer--Cartan forms.
In Section~\ref{moving-frames-sec}, we begin the method of moving frames and identify f\/irst-order invariants for
centroaf\/f\/ine surfaces.
Based on these invariants, nondegenerate surfaces may be locally classif\/ied as ``space-like'', ``time-like'', or ``null''.
In Sections~\ref{space-like-sec} and~\ref{time-like-sec}, we continue the method of moving frames for the space-like and
time-like cases, respectively.
(We do not consider the null case here; it may be explored in a~future paper.) Finally, in Section~\ref{homog-sec} we
classify the homogeneous examples in both the space-like and time-like cases.

\section[Centroaf\/f\/ine surfaces in $\mathbb{R}^5$, adapted frames, and Maurer--Cartan forms]{Centroaf\/f\/ine surfaces
in $\boldsymbol{\mathbb{R}^5}$, adapted frames,\\ and Maurer--Cartan forms}
\label{background-sec}

Five-dimensional {\em centroaffine space} is the manifold $\mathbb{R}^5 \setminus \{0\}$, equipped with a~natural ${\rm
GL}(5,\mathbb{R})$-action.
Specif\/ically, ${\rm GL}(5, \mathbb{R})$ acts on $\mathbb{R}^5 \setminus \{0\}$ by left multiplication: for $\mathbf{x}
\in \mathbb{R}^5 \setminus \{0\}$, $g \in {\rm GL}(5, \mathbb{R})$, we have
\begin{gather*}
g \cdot \mathbf{x} = g\mathbf{x}.
\end{gather*}
The group ${\rm GL}(5,\mathbb{R})$ may be regarded as a~principal bundle over $\mathbb{R}^5 \setminus \{0\}$: write an
arbitrary element $g\in {\rm GL}(5,\mathbb{R})$ as
\begin{gather*}
g =
\begin{bmatrix}
\mathbf{e}_0 & \mathbf{e}_1 & \mathbf{e}_2 & \mathbf{e}_3 & \mathbf{e}_4
\end{bmatrix},
\end{gather*}
where $\mathbf{e}_0, \ldots, \mathbf{e}_4 \in \mathbb{R}^5 \setminus \{0\}$ are linearly independent column vectors.
Then def\/ine the bundle map $\pi:{\rm GL}(5,\mathbb{R}) \to \mathbb{R}^5 \setminus \{0\}$~by
\begin{gather}
\pi \left(
\begin{bmatrix}
\mathbf{e}_0 & \mathbf{e}_1 & \mathbf{e}_2 & \mathbf{e}_3 & \mathbf{e}_4
\end{bmatrix}
\right) = \mathbf{e}_0.
\label{define-pi}
\end{gather}
The f\/iber group~$H$ is isomorphic to the stabilizer of the point
\begin{gather*}
\begin{bmatrix}
1 & 0 & 0 & 0 & 0
\end{bmatrix}^{\rm T} \in \mathbb{R}^5 \setminus \{0\},
\end{gather*}
and this construction endows the manifold $\mathbb{R}^5 \setminus \{0\}$ with the structure of the homogeneous space
${\rm GL}(5,\mathbb{R})/H$.

We also think of the bundle $\pi:{\rm GL}(5,\mathbb{R}) \to \mathbb{R}^5 \setminus \{0\}$ as the {\em centroaffine frame
bundle} $\mathcal{F}$ over $\mathbb{R}^5 \setminus \{0\}$.
For each point $\mathbf{e}_0 \in \mathbb{R}^5 \setminus \{0\}$, the f\/iber over $\mathbf{e}_0$ consists of all frames
$(\mathbf{e}_0, \mathbf{e}_1, \mathbf{e}_2, \mathbf{e}_3, \mathbf{e}_4)$ for the tangent space
$T_{\mathbf{e}_0}(\mathbb{R}^5 \setminus \{0\})$~-- i.e., all frames for which the f\/irst vector in the frame is equal to
the position vector.

The {\em Maurer--Cartan forms} $\omega^i_j$ on $\mathcal{F}$ are def\/ined by the equations
\begin{gather}
d\mathbf{e}_{i} = \mathbf{e}_{j}  \omega^{j}_{i},
\qquad
0 \leq i,j \leq 4,
\label{define-MC-forms}
\end{gather}
and they satisfy the {\em Cartan structure equations}
\begin{gather}
d\omega^i_j = -\omega^i_k \wedge \omega^k_j.
\label{structure-eqns}
\end{gather}
(For details, see~\cite[p.~18]{CFB} or~\cite{FTC}.)

We are interested in the geometry of 2-dimensional immersions $\bar{f}:M^2 \to \mathbb{R}^5 \setminus \{0\}$; we will
use Cartan's method of moving frames to compute local invariants for such immersions under the action of ${\rm
GL}(5,\mathbb{R})$.

\begin{dfn}
\label{centroaffine-immersion-def}
An immersion $\bar{f}$ of a~2-dimensional manifold~$M$ into $\mathbb{R}^5 \setminus \{0\}$ is called a~{\em centroaffine
immersion} if the position vector $\bar{f}(x)$ is transversal to the tangent space $\bar{f}_*(T_xM)$ for all $x\in M$.
The image $\Sigma = \bar{f}(M)$ is called a~{\em centroaffine surface} in $\mathbb{R}^5 \setminus \{0\}$.
\end{dfn}

In order to begin the method of moving frames, consider the induced bundle of centroaf\/f\/ine frames along $\Sigma =
\bar{f}(M)$; this is simply the pullback bundle $\mathcal{F}_0 = \bar{f}^*\mathcal{F}$ over~$M$.
A~{\em centroaffine frame field} along~$\Sigma$ is a~section of $\mathcal{F}_0$~-- i.e., a~smooth map $f:M \to {\rm
GL}(5,\mathbb{R})$ such that $\pi \circ f = \bar{f}$. Throughout the remainder of this paper, we will consider the
pullbacks of the Maurer--Cartan forms on $\mathcal{F}$ to~$M$ via such sections~$f$, and we will suppress the pullback
notation.

We will gradually adapt our choice of centroaf\/f\/ine frame f\/ields based on the geometry of~$\Sigma$.
For our f\/irst adaptation, consider the subbundle $\mathcal{F}_1 \subset \mathcal{F}_0$ consisting of all frames for
which ($\mathbf{e}_1(x)$, $\mathbf{e}_2(x)$) span the tangent space $T_{\bar{f}(x)}\Sigma$ for each $x \in M$.
A~section $f:M \to \mathcal{F}_1$ will be called a~{\em $1$-adapted frame field} along~$\Sigma$.
Any two 1-adapted frames $(\mathbf{e}_0, \ldots, \mathbf{e}_4)$, $(\tilde{\mathbf{e}}_0, \ldots, \tilde{\mathbf{e}}_4)$
based at the same point $x \in M$ are related by a~transformation of the form
\begin{gather}
\begin{bmatrix}
\tilde{\mathbf{e}}_0 & \tilde{\mathbf{e}}_1 & \tilde{\mathbf{e}}_2 & \tilde{\mathbf{e}}_3 & \tilde{\mathbf{e}}_4
\end{bmatrix}
=
\begin{bmatrix}
\mathbf{e}_0 & \mathbf{e}_1 & \mathbf{e}_2 & \mathbf{e}_3 & \mathbf{e}_4
\end{bmatrix}
\begin{bmatrix}
1 & 0 & 0 & r_{03} & r_{04}
\\
0 & a_{11} & a_{12} & r_{13} & r_{14}
\\
0 & a_{21} & a_{22} & r_{23} & r_{24}
\\
0 & 0 & 0 & b_{33} & b_{34}
\\
0 & 0 & 0 & b_{43} & b_{44}
\end{bmatrix},
\label{1-adapted-frame-transformation}
\end{gather}
where the $2\times 2$ submatrices
\begin{gather*}
A =
\begin{bmatrix}
a_{11} & a_{12}
\\
a_{21} & a_{22}
\end{bmatrix},
\qquad
B =
\begin{bmatrix}
b_{33} & b_{34}
\\
b_{43} & b_{44}
\end{bmatrix}
\end{gather*}
are elements of ${\rm GL}(2, \mathbb{R})$.
We will denote the group of all matrices of the form in~\eqref{1-adapted-frame-transformation} by $G_1$; then the bundle
$\mathcal{F}_1$ is a~principal bundle over~$M$ with f\/iber group $G_1$.

If $f, \tilde{f}:M \to \mathcal{F}_1$ are two 1-adapted frame f\/ields along~$\Sigma$, then $f$, $\tilde{f}$ are related~by
the equation
\begin{gather*}
\tilde{f}(x) = f(x) \cdot g(x)
\end{gather*}
for some smooth function $g:M \to G_1$, as in equation~\eqref{1-adapted-frame-transformation}.
Then the corresponding $\mathfrak{gl}(5, \mathbb{R})$-valued Maurer--Cartan forms $\Omega = [\omega^i_j]$,
$\tilde{\Omega} = [\tilde{\omega}^i_j]$ on~$M$ are related as follows:
\begin{gather}
\tilde{\Omega} = g^{-1} dg + g^{-1} \Omega g.
\label{how-MC-forms-transform}
\end{gather}

\section{Reduction of the structure group and f\/irst-order invariants}\label{moving-frames-sec}

Now consider the pullbacks of the Maurer--Cartan forms to~$M$ via a~1-adapted frame f\/ield~$f$.
From equation~\eqref{define-MC-forms} for $d\mathbf{e}_0$ and the fact that the image of $d\mathbf{e}_0$ is spanned~by
$\mathbf{e}_1$ and $\mathbf{e}_2$, we have
\begin{gather}
\omega^0_0 = \omega^3_0 = \omega^4_0 = 0.
\label{first-vanishing-conds}
\end{gather}
Moreover, the 1-forms $\omega^1_0$, $\omega^2_0$ are semi-basic for the projection $\pi:\mathcal{F}_1 \to M$; in fact,
they form a~basis for the semi-basic 1-forms on $\mathcal{F}_1$.

Dif\/ferentiating equations~\eqref{first-vanishing-conds} yields:
\begin{gather*}
0 = d\omega^0_0 = -\big(\omega^0_1 \wedge \omega^1_0 + \omega^0_2 \wedge \omega^2_0\big),
\\
0 = d\omega^3_0 = -\big(\omega^3_1 \wedge \omega^1_0 + \omega^3_2 \wedge \omega^2_0\big),
\\
0 = d\omega^4_0 = -\big(\omega^4_1 \wedge \omega^1_0 + \omega^4_2 \wedge \omega^2_0\big).
\end{gather*}
Cartan's lemma (see, e.g.,~\cite{CFB} or~\cite{FTC}) then implies that there exist functions $h^k_{ij} = h^k_{ji}$, $k=0,3,4$, on~$M$ such that
\begin{gather*}
\begin{bmatrix}
\omega^k_1
\vspace{1mm}\\
\omega^k_2
\end{bmatrix}
=
\begin{bmatrix}
h^k_{11} & h^k_{12}
\vspace{1mm}\\
h^k_{12} & h^k_{22}
\end{bmatrix}
\begin{bmatrix}
\omega^1_0
\vspace{1mm}\\
\omega^2_0
\end{bmatrix}.
\end{gather*}
For simplicity of notation, let $h^k$ denote the matrix
\begin{gather*}
h^k =
\begin{bmatrix}
h^k_{11} & h^k_{12}
\vspace{1mm}\\
h^k_{12} & h^k_{22}
\end{bmatrix},
\qquad
k=0,3,4.
\end{gather*}

If $f, \tilde{f}:M \to \mathcal{F}_1$ are two 1-adapted frame f\/ields related by a~transformation of the
form~\eqref{1-adapted-frame-transformation}, then we can use equation~\eqref{how-MC-forms-transform} to determine how
the corresponding matrices $h^k$, $\tilde{h}^k$ are related.
First, it follows from the fact that $\tilde{\mathbf{e}}_0 = \mathbf{e}_0$ that
\begin{gather*}
d\tilde{\mathbf{e}}_0 =
\begin{bmatrix}
\tilde{\mathbf{e}}_1 & \tilde{\mathbf{e}}_2
\end{bmatrix}
\begin{bmatrix}
\tilde{\omega}^1_0
\vspace{1mm}\\
\tilde{\omega}^2_0
\end{bmatrix}
=
\begin{bmatrix}
\mathbf{e}_1 & \mathbf{e}_2
\end{bmatrix}
\begin{bmatrix}
\omega^1_0
\vspace{1mm}\\
\omega^2_0
\end{bmatrix}
= d\mathbf{e}_0.
\end{gather*}
Then, since we have
\begin{gather*}
\begin{bmatrix}
\tilde{\mathbf{e}}_1 & \tilde{\mathbf{e}}_2
\end{bmatrix}
=
\begin{bmatrix}
\mathbf{e}_1 & \mathbf{e}_2
\end{bmatrix}
A,
\end{gather*}
we must have
\begin{gather}
\begin{bmatrix}
\tilde{\omega}^1_0
\vspace{1mm}\\
\tilde{\omega}^2_0
\end{bmatrix}
= A^{-1}
\begin{bmatrix}
\omega^1_0
\vspace{1mm}\\
\omega^2_0
\end{bmatrix}.
\label{semi-basic-trans}
\end{gather}
Similar considerations show that
\begin{gather}
\begin{bmatrix}
\tilde{\omega}^3_1
\vspace{1mm}\\
\tilde{\omega}^3_2
\end{bmatrix}
= \frac{1}{(\det B)} A^{\rm T} \left(b_{44}
\begin{bmatrix}
\omega^3_1
\vspace{1mm}\\
\omega^3_2
\end{bmatrix}
- b_{34}
\begin{bmatrix}
\omega^4_1
\vspace{1mm}\\
\omega^4_2
\end{bmatrix}
\right),
\nonumber
\\
\begin{bmatrix}
\tilde{\omega}^4_1
\vspace{1mm}\\
\tilde{\omega}^4_2
\end{bmatrix}
= \frac{1}{(\det B)} A^{\rm T} \left(-b_{43}
\begin{bmatrix}
\omega^3_1
\vspace{1mm}\\
\omega^3_2
\end{bmatrix}
+ b_{33}
\begin{bmatrix}
\omega^4_1
\vspace{1mm}\\
\omega^4_2
\end{bmatrix}
\right),
\nonumber
\\
\begin{bmatrix}
\tilde{\omega}^0_1
\vspace{1mm}\\
\tilde{\omega}^0_2
\end{bmatrix}
= A^{\rm T}
\begin{bmatrix}
\omega^0_1
\vspace{1mm}\\
\omega^0_2
\end{bmatrix}
- r_{03}
\begin{bmatrix}
\tilde{\omega}^3_1
\vspace{1mm}\\
\tilde{\omega}^3_2
\end{bmatrix}
- r_{04}
\begin{bmatrix}
\tilde{\omega}^4_1
\vspace{1mm}\\
\tilde{\omega}^4_2
\end{bmatrix}.
\label{1-adapted-connection-trans}
\end{gather}
Together, equations~\eqref{semi-basic-trans},~\eqref{1-adapted-connection-trans} imply that
\begin{gather}
\tilde{h}^3 = \frac{1}{(\det B)} A^{\rm T} \left(b_{44} h^3 - b_{34} h^4 \right) A,
\qquad
\tilde{h}^4 = \frac{1}{(\det B)} A^{\rm T} \left(-b_{43} h^3 + b_{33} h^4 \right) A,
\nonumber
\\
\tilde{h}^0 = A^{\rm T} h^0 A-r_{03} \tilde{h}^3 - r_{04} \tilde{h}^4.
\label{1-adapted-group-action}
\end{gather}

\begin{dfn}
\label{nondegenerate-def}
A~centroaf\/f\/ine surface $\Sigma = \bar{f}(M)$ will be called {\em nondegenerate} if the matri\-ces~$h^0$, $h^3$, $h^4$ are
linearly independent in ${\rm Sym}^2(\mathbb{R})$ at every point of~$M$.
\end{dfn}

Henceforth, we assume that~$\Sigma$ is nondegenerate; from the group action~\eqref{1-adapted-group-action} it is clear
that this def\/inition is independent of the choice of 1-adapted frame f\/ield $f:M \to \mathcal{F}_1$ along~$\Sigma$.

The next step is to use the group action~\eqref{1-adapted-group-action} to f\/ind normal forms for the matri\-ces~$h^0$,~$h^3$,~$h^4$.
First consider the action on~$h^3$,~$h^4$: it can be written as the composition of two separate actions by the matrices
$A,B \in {\rm GL}(2, \mathbb{R})$:
\begin{gather*}
A\cdot \big(h^3, h^4\big)=\big(A^{\rm T} h^3 A, A^{\rm T} h^4 A\big),
\\
B \cdot \big(h^3, h^4\big)=\left(\frac{1}{(\det B)} \left(b_{44} h^3 - b_{34} h^4 \right), \frac{1}{(\det
B)}\left(-b_{43} h^3 + b_{33} h^4 \right)\right).
\end{gather*}
If we let~$P$ denote the 2-dimensional subspace of ${\rm Sym}^2(\mathbb{R})$ spanned by $(h^3, h^4)$, then we see that
the action by~$B$ preserves~$P$, while~$A$ acts on~$P$ via
\begin{gather}
A\cdot P = A^{\rm T} P A.
\label{plane-action}
\end{gather}

In order to understand the action~\eqref{plane-action}, consider the related action
\begin{gather}
A\cdot h = A^{\rm T} h A
\label{sym-2-action}
\end{gather}
on ${\rm Sym}^2(\mathbb{R})$.
It is shown in~\cite[p.~115]{MH73} that this action preserves the indef\/inite quadratic form
\begin{gather}
Q(h) = -\det(h)
\label{quadratic-form}
\end{gather}
up to a~scale factor.
More precisely,~$Q$ gives ${\rm Sym}^2(\mathbb{R})$ the structure of the Minkowski space $\mathbb{R}^{2,1}$, and the
action~\eqref{sym-2-action} gives a~representation of ${\rm GL}(2,\mathbb{R})$ as the group ${\rm CSO}^+(2,1)$ of
orientation-preserving, orthochronous, conformal Minkowski transformations of $\mathbb{R}^{2,1}$.
This action has precisely 6 orbits, represented by the matrices
\begin{gather}
\label{matrix-reps}
\begin{bmatrix}
0 & 0
\\
0 & 0
\end{bmatrix},
\qquad
\begin{bmatrix}
1 & 0
\\
0 & 0
\end{bmatrix},
\qquad
\begin{bmatrix}
-1 & 0
\\
0 & 0
\end{bmatrix},
\qquad
\begin{bmatrix}
1 & 0
\\
0 & 1
\end{bmatrix},
\qquad
\begin{bmatrix}
1 & 0
\\
0 & -1
\end{bmatrix},
\qquad
\begin{bmatrix}
-1 & 0
\\
0 & -1
\end{bmatrix}.
\end{gather}
In terms of the Minkowski metric on ${\rm Sym}^2(\mathbb{R})$ determined by~$Q$, the second and third matrices
in~\eqref{matrix-reps} represent the two orbits of null vectors (oriented in opposite time directions); the fourth and
f\/ifth represent the two orbits of time-like vectors (oriented in opposite time directions), while the sixth represents
the single orbit of space-like vectors.

The action~\eqref{sym-2-action} induces the action~\eqref{plane-action} on the Grassmannian of 2-planes in ${\rm Sym}^2(\mathbb{R})$,
which corresponds to the action of ${\rm CSO}^+(2,1)$ on planes in $\mathbb{R}^{2,1}$.
This action has precisely 3~orbits, consisting of planes that are space-like, time-like, or null with respect to the
quadratic form~\eqref{quadratic-form}.
These orbits are represented by the 2-planes
\begin{gather}
P_1   = \text{span} \left( \begin{bmatrix} 1 & 0 \\  0 & -1 \end{bmatrix}, \begin{bmatrix} 0 & 1 \\  1 & 0 \end{bmatrix} \right)   \qquad   \text{(space-like)}, \nonumber\\
P_2   = \text{span} \left( \begin{bmatrix} 1 & 0 \\  0 & 1 \end{bmatrix}, \begin{bmatrix} 1 & 0 \\  0 & -1 \end{bmatrix} \right)   \qquad   \text{(time-like)}, \nonumber\\
P_3   = \text{span} \left( \begin{bmatrix} 1 & 0 \\  0 & 0 \end{bmatrix}, \begin{bmatrix} 0 & 1 \\  1 & 0 \end{bmatrix} \right)   \qquad    \text{(null)}.
 \label{normal-forms-planes}
\end{gather}
The {\em type} of the plane~$P$ spanned by $(h^3, h^4)$ (space-like, time-like, or null) is preserved by the group
action~\eqref{1-adapted-group-action}; thus the type of~$P$ at any point $x \in M$ is well-def\/ined, independent of the
choice of 1-adapted frame f\/ield $f:M \to \mathcal{F}_1$.

\begin{rmk}
We note that similar quadratic forms $h_3$, $h_4$ were obtained by Nomizu and Vrancken in their study of af\/f\/ine surfaces
in $\mathbb{R}^4$; see~\cite{NV93}.
\end{rmk}

At this point, the method of moving frames dictates that we divide into cases based on the type of~$P$.
In order to proceed, we make the following assumption:

\begin{assumption}
\label{constant-type-assumption}
Assume that~$\Sigma$ has {\em constant type}~-- i.e., that the type of~$P$ is the same at every point $x \in M$.
\end{assumption}

In this paper we will consider only the space-like and time-like cases; the null case is considerably more complicated and
may be explored in a~future paper.

\section{The space-like case}\label{space-like-sec}

First, suppose that the plane~$P$ spanned by $(h^3, h^4)$ is space-like at every point of~$M$.
According to the group action~\eqref{plane-action}, we can f\/ind a~1-adapted frame f\/ield along~$\Sigma$ for which
$P=P_1$, as in equation~\eqref{normal-forms-planes}.
Furthermore, we can then use the action by~$B$ to f\/ind a~1-adapted frame f\/ield along~$\Sigma$ for which
\begin{gather}
\label{space-like-h3-h4}
h^3 =
\begin{bmatrix}
1 & 0
\\
0 & -1
\end{bmatrix},
\qquad
h^4 =
\begin{bmatrix}
0 & 1
\\
1 & 0
\end{bmatrix}.
\end{gather}

The next step is to determine the subgroup of $G_1$ that preserves the conditions~\eqref{space-like-h3-h4}.
To this end, f\/irst note that the plane $P_1$ spanned by the matrices~\eqref{space-like-h3-h4} consists precisely of the
trace-free matrices in ${\rm Sym}^2(\mathbb{R})$.
A~straightforward computation shows that, with $h^3$, $h^4$ as in~\eqref{space-like-h3-h4},
\begin{gather*}
\tr\big(A^{\rm T} h^3 A\big) = a_{11}^2 + a_{12}^2 - a_{21}^2 - a_{22}^2,
\qquad
\tr\big(A^{\rm T} h^4 A\big) = 2(a_{11} a_{21} + a_{12} a_{22}).
\end{gather*}
Therefore, the action~\eqref{plane-action} preserves $P_1$ if and only if
\begin{gather*}
a_{11}^2 + a_{12}^2 - a_{21}^2 - a_{22}^2 = a_{11} a_{21} + a_{12} a_{22} = 0,
\end{gather*}
which is true if and only if
\begin{gather*}
A = \lambda A_0
\end{gather*}
for some $\lambda \in \mathbb{R}^*$, $A_0 \in {\rm O}(2, \mathbb{R})$.
Substituting this condition into equations~\eqref{1-adapted-group-action} and imposing the conditions $\tilde{h}^3 =
h^3$, $\tilde{h}^4 = h^4$ then yields
\begin{gather*}
B = \lambda^2 (A_0)^2.
\end{gather*}

For simplicity, we will restrict to transformations with $A_0 \in {\rm SO}(2,\mathbb{R})$, $\lambda > 0$; this has the
advantage of producing a~frame bundle whose f\/iber is a~connected Lie group.
Thus we will assume that
\begin{gather}
A=
\begin{bmatrix}
\lambda \cos (\theta) & -\lambda \sin (\theta)
\\
\lambda \sin (\theta) & \lambda \cos (\theta)
\end{bmatrix},
\qquad
B =
\begin{bmatrix}
\lambda^2 \cos (2\theta) & -\lambda^2 \sin (2\theta)
\\
\lambda^2 \sin (2\theta) & \lambda^2 \cos (2\theta)
\end{bmatrix},
\label{space-like-AB}
\end{gather}
where $\lambda > 0$, $\theta \in \mathbb{R}$.

Next, consider the ef\/fect of the action~\eqref{1-adapted-group-action} on $h^0$.
With $A = I_2$ and $r_{03}$, $r_{04}$ chosen appropriately, we can add any linear combination of $h^3$, $h^4$ to $h^0$.
Thus we can f\/ind a~1-adapted frame f\/ield for which $h^0$ is a~multiple (nonzero by the nondegeneracy assumption) of the
identity matrix $I_2$.
Then under the action~\eqref{1-adapted-group-action} with $A$, $B$ as in~\eqref{space-like-AB}, we have
\begin{gather*}
\tilde{h}^0 = \lambda^2 h^0 +
\begin{bmatrix}
-r_{03} & -r_{04}
\\
-r_{04} & r_{03}
\end{bmatrix}.
\end{gather*}
Therefore, we can f\/ind a~1-adapted frame f\/ield along~$\Sigma$ satisfying the additional condition that $h^0 = \pm I_2$,
and this condition is preserved by transformations of the form~\eqref{1-adapted-group-action} with $A$, $B$ as
in~\eqref{space-like-AB}, $\lambda=1$, and and $r_{03} = r_{04} = 0$.

\begin{dfn}
Let $\Sigma = \bar{f}(M)$ be a~nondegenerate, space-like centroaf\/f\/ine surface in $\mathbb{R}^5 \setminus \{0\}$.
A~1-adapted frame f\/ield $f:M \to \mathcal{F}_1$ will be called {\em $2$-adapted} if it satisf\/ies the conditions
\begin{gather}
h^3 =
\begin{bmatrix}
1 & 0
\\
0 & -1
\end{bmatrix},
\qquad
h^4 =
\begin{bmatrix}
0 & 1
\\
1 & 0
\end{bmatrix},
\qquad
h^0 =
\begin{bmatrix}
\epsilon & 0
\\
0 & \epsilon
\end{bmatrix},
\label{2-adapted-space-like-h}
\end{gather}
with $\epsilon=\pm 1$, at every point of~$M$.
\end{dfn}

Any two 2-adapted frames $(\mathbf{e}_0, \ldots, \mathbf{e}_4)$, $(\tilde{\mathbf{e}}_0, \ldots, \tilde{\mathbf{e}}_4)$
based at the same point $x \in M$ are related by a~transformation of the form
\begin{gather}
\begin{bmatrix}
\tilde{\mathbf{e}}_0 & \tilde{\mathbf{e}}_1 & \tilde{\mathbf{e}}_2 & \tilde{\mathbf{e}}_3 & \tilde{\mathbf{e}}_4
\end{bmatrix}
\nonumber
\\
\qquad{}
=
\begin{bmatrix}
\mathbf{e}_0 & \mathbf{e}_1 & \mathbf{e}_2 & \mathbf{e}_3 & \mathbf{e}_4
\end{bmatrix}
\begin{bmatrix}
1 & 0 & 0 & 0 & 0
\\
0 & \cos(\theta) & -\sin(\theta) & r_{13} & r_{14}
\\
0 & \sin(\theta) & \cos(\theta) & r_{23} & r_{24}
\\
0 & 0 & 0 & \cos(2\theta) & -\sin(2\theta)
\\
0 & 0 & 0 & \sin(2\theta) & \cos(2\theta)
\end{bmatrix}.
\label{2-adapted-space-like-frame-transformation}
\end{gather}
We will denote the group of all matrices of the form in~\eqref{2-adapted-space-like-frame-transformation} by $G_2$; then
the 2-adapted frame f\/ields along~$\Sigma$ are the smooth sections of a~principal bundle~$\mathcal{F}_2 \subset
\mathcal{F}_1$ over~$M$ with f\/iber group~$G_2$.

The equations~\eqref{2-adapted-space-like-h} are equivalent to the condition that the Maurer--Cartan forms associated to
a~2-adapted frame f\/ield satisfy the conditions
\begin{gather}
\omega^3_1= \omega^1_0,
\qquad
\omega^3_2=-\omega^2_0,
\qquad
\omega^4_1=\omega^2_0,
\qquad
\omega^4_2=\omega^1_0,
\qquad
\omega^0_1=\epsilon  \omega^1_0,
\qquad
\omega^0_2=\epsilon  \omega^2_0.
\label{2-adapted-space-like-MC-relations}
\end{gather}
Dif\/ferentiating equations~\eqref{2-adapted-space-like-MC-relations} yields
\begin{gather*}
\big(2\omega^1_1 - \omega^3_3\big) \wedge \omega^1_0 + \big(\omega^1_2 - \omega^2_1 - \omega^3_4\big) \wedge \omega^2_0 = 0,
\\
\big(\omega^1_2 - \omega^2_1 - \omega^3_4\big) \wedge \omega^1_0 + \big(\omega^3_3 - 2\omega^2_2\big) \wedge \omega^2_0 = 0,
\\
\big(2\omega^2_1 - \omega^4_3\big) \wedge \omega^1_0 + \big(\omega^1_1 + \omega^2_2 - \omega^4_4\big) \wedge \omega^2_0 = 0,
\\
\big(\omega^1_1 + \omega^2_2 - \omega^4_4\big) \wedge \omega^1_0 + \big(2\omega^1_2 + \omega^4_3\big) \wedge \omega^2_0 = 0,
\\
\big(2\epsilon\omega^1_1 - \omega^0_3\big) \wedge \omega^1_0 + \big(\epsilon\omega^1_2 + \epsilon\omega^2_1 - \omega^0_4\big)
\wedge \omega^2_0 = 0,
\\
\big(\epsilon\omega^1_2 + \epsilon\omega^2_1 - \omega^0_4\big) \wedge \omega^1_0 + \big(2\epsilon\omega^2_2 + \omega^0_3\big)
\wedge \omega^2_0 = 0.
\end{gather*}
Applying Cartan's lemma to these equations shows that there exists a~1-form~$\alpha$ and func\-tions~$h^i_{jk}$ on
$\mathcal{F}_2$ such that
\begin{alignat}{3}
& \omega^0_3=h^0_{31}\omega^1_0 + h^0_{32}\omega^2_0,
\qquad&&
\omega^0_4=h^0_{41}\omega^1_0 + h^0_{42}\omega^2_0,&
\nonumber
\\
&\omega^1_1=h^1_{11}\omega^1_0 + h^1_{12}\omega^2_0,
\qquad&&
\omega^1_2=\alpha + h^1_{21}\omega^1_0 + h^1_{22}\omega^2_0,&
\nonumber
\\
&\omega^2_1=-\alpha + h^1_{21}\omega^1_0 + h^1_{22}\omega^2_0,
\qquad&&
\omega^2_2=h^2_{21}\omega^1_0 + h^2_{22}\omega^2_0,&
\nonumber
\\
&\omega^3_3=h^3_{31}\omega^1_0 + h^3_{32}\omega^2_0,
\qquad&&
\omega^3_4=2\alpha + h^3_{41}\omega^1_0 + h^3_{42}\omega^2_0,&
\nonumber
\\
&\omega^4_3=-2\alpha + h^4_{31}\omega^1_0 + h^4_{32}\omega^2_0,
\qquad&&
\omega^4_4=h^4_{41}\omega^1_0 + h^4_{42}\omega^2_0.&
\label{2-adapted-space-like-MC-semi-basic-forms}
\end{alignat}
Moreover, the functions $h^i_{jk}$ satisfy the relations
\begin{alignat}{3}
& 2 h^1_{12} - h^3_{32} + h^3_{41} = 0,
\qquad &&
2 h^2_{21} - h^3_{31} - h^3_{42} = 0,&
\nonumber
\\
& h^1_{11} - 2 h^1_{22} + h^2_{21} + h^4_{32} - h^4_{41} = 0,
\qquad &&
h^1_{12} - 2 h^1_{21} + h^2_{22} - h^4_{31} - h^4_{42} = 0,&
\nonumber
\\
& h^0_{32} - h^0_{41} + 2\epsilon\big(h^1_{21} - h^1_{12}\big) = 0,
\qquad &&
h^0_{31} + h^0_{42} + 2\epsilon\big(h^2_{21} - h^1_{22}\big) = 0.&
\label{2-adapted-space-like-h-relations}
\end{alignat}

If $f, \tilde{f}:M \to \mathcal{F}_2$ are two 2-adapted frame f\/ields related by a~transformation of the
form~\eqref{2-adapted-space-like-frame-transformation}, then we can once again use
equation~\eqref{how-MC-forms-transform} to determine how the corresponding functions $h^i_{jk}$, $\tilde{h}^i_{jk}$ are
related.
Some of these relationships are more complicated than others; the most straightforward to compute are those
corresponding to the forms $\tilde{\omega}^0_3$, $\tilde{\omega}^0_4$.
These forms appear as the coef\/f\/icients of $\tilde{\mathbf{e}}_0 = \mathbf{e}_0$ in the equations~\eqref{define-MC-forms}
for $d\tilde{\mathbf{e}}_3$, $d\tilde{\mathbf{e}}_4$.
By applying equations~\eqref{2-adapted-space-like-frame-transformation}
and~\eqref{2-adapted-space-like-MC-semi-basic-forms}, one can show that
\begin{gather*}
\begin{bmatrix}
\tilde{h}^0_{31}
\vspace{1mm}\\
\tilde{h}^0_{41}
\end{bmatrix}
=
\begin{bmatrix}
\cos(2\theta) & \sin(2\theta)
\\
-\sin(2\theta) & \cos(2\theta)
\end{bmatrix}
\begin{bmatrix}
h^0_{31}
\vspace{1mm}\\
h^0_{41}
\end{bmatrix}
+ \epsilon
\begin{bmatrix}
r_{13}
\\
r_{14}
\end{bmatrix},
\\
\begin{bmatrix}
\tilde{h}^0_{32}
\vspace{1mm}\\
\tilde{h}^0_{42}
\end{bmatrix}
=
\begin{bmatrix}
\cos(2\theta) & \sin(2\theta)
\\
-\sin(2\theta) & \cos(2\theta)
\end{bmatrix}
\begin{bmatrix}
h^0_{32}
\vspace{1mm}\\
h^0_{42}
\end{bmatrix}
+ \epsilon
\begin{bmatrix}
r_{23}
\\
r_{24}
\end{bmatrix}.
\end{gather*}
Thus we can f\/ind a~2-adapted frame f\/ield along~$\Sigma$ satisfying the conditions that
\begin{gather*}
h^0_{31} = h^0_{32} = h^0_{41} = h^0_{42} = 0,
\end{gather*}
and these conditions are preserved by transformations of the form~\eqref{2-adapted-space-like-frame-transformation} with
$r_{13} = r_{14} = r_{23} = r_{24} = 0$.

\begin{dfn}
Let $\Sigma = \bar{f}(M)$ be a~nondegenerate, space-like centroaf\/f\/ine surface in $\mathbb{R}^5 \setminus \{0\}$.
A~2-adapted frame f\/ield $f:M \to \mathcal{F}_2$ will be called {\em $3$-adapted} if it satisf\/ies the conditions
\begin{gather}
h^0_{31} = h^0_{32} = h^0_{41} = h^0_{42} = 0
\label{3-adapted-space-like-h}
\end{gather}
at every point of~$M$.
\end{dfn}

Any two 3-adapted frames $(\mathbf{e}_0, \ldots, \mathbf{e}_4)$, $(\tilde{\mathbf{e}}_0, \ldots, \tilde{\mathbf{e}}_4)$
based at the same point $x \in M$ are related by a~transformation of the form
\begin{gather}
\begin{bmatrix}
\tilde{\mathbf{e}}_0 & \tilde{\mathbf{e}}_1 & \tilde{\mathbf{e}}_2 & \tilde{\mathbf{e}}_3 & \tilde{\mathbf{e}}_4
\end{bmatrix}
\nonumber
\\
\qquad{}
=
\begin{bmatrix}
\mathbf{e}_0 & \mathbf{e}_1 & \mathbf{e}_2 & \mathbf{e}_3 & \mathbf{e}_4
\end{bmatrix}
\begin{bmatrix}
1 & 0 & 0 & 0 & 0
\\
0 & \cos(\theta) & -\sin(\theta) & 0 & 0
\\
0 & \sin(\theta) & \cos(\theta) & 0 & 0
\\
0 & 0 & 0 & \cos(2\theta) & -\sin(2\theta)
\\
0 & 0 & 0 & \sin(2\theta) & \cos(2\theta)
\end{bmatrix}.
\label{3-adapted-space-like-frame-transformation}
\end{gather}
We will denote the group of all matrices of the form in~\eqref{3-adapted-space-like-frame-transformation} by $G_3$; note
that $G_3 \cong {\rm SO}(2,\mathbb{R})$.
Then the 3-adapted frame f\/ields are the smooth sections of a~principal bundle $\mathcal{F}_3 \subset \mathcal{F}_2$
over~$M$ with f\/iber group $G_3$.

The equations~\eqref{3-adapted-space-like-h} are equivalent to the condition that the Maurer--Cartan forms associated to
a~3-adapted frame f\/ield satisfy the conditions
\begin{gather}
\omega^0_3 = \omega^0_4 = 0.
\label{3-adapted-MC-space-like-relations}
\end{gather}
Dif\/ferentiating equations~\eqref{3-adapted-MC-space-like-relations} yields
\begin{gather*}
\omega^1_3 \wedge \omega^1_0 + \omega^2_3 \wedge \omega^2_0 = \omega^1_4 \wedge \omega^1_0 + \omega^2_4 \wedge
\omega^2_0 = 0,
\end{gather*}
and applying Cartan's lemma shows that there exist functions $h^i_{jk}$ on $\mathcal{F}_3$ such that
\begin{alignat*}{3}
& \omega^1_3=h^1_{31}\omega^1_0 + h^1_{32}\omega^2_0,
\qquad&&
\omega^2_3=h^1_{32}\omega^1_0 + h^2_{32}\omega^2_0,&
\\
&\omega^1_4=h^1_{41}\omega^1_0 + h^1_{42}\omega^2_0,
\qquad&&
\omega^2_4=h^1_{42}\omega^1_0 + h^2_{42}\omega^2_0.&
\end{alignat*}
Moreover, on $\mathcal{F}_3$, the last two relations in equations~\eqref{2-adapted-space-like-h-relations} simplify to
\begin{gather*}
h^1_{21} - h^1_{12} = h^2_{21} - h^1_{22} = 0.
\end{gather*}

At this point, we have canonically associated to any nondegenerate, space-like centroaf\/f\/ine surface in $\mathbb{R}^5
\setminus \{0\}$ a~frame bundle $\mathcal{F}_3$ over~$M$ with f\/iber group isomorphic to ${\rm SO}(2,\mathbb{R})$.
Thus we have the following theorem:

\begin{teo}
\label{space-like-thm}
Let $\bar{f}:M \to \mathbb{R}^5 \setminus\{0\}$ be a~centroaffine immersion whose image $\Sigma = \bar{f}(M)$ is
a~nondegenerate, space-like centroaffine surface.
Then the pullbacks of the Maurer--Cartan forms on ${\rm GL}(5,\mathbb{R})$ to the bundle $\mathcal{F}_3$ of $3$-adapted
frames on~$\Sigma$ determine a~well-defined Riemannian metric
\begin{gather*}
I = \left(\omega^1_0\right)^2 + \left(\omega^2_0\right)^2
\end{gather*}
on~$\Sigma$, called the centroaffine metric.
Moreover, there is a~well-defined ``centroaffine normal bundle'' $N\Sigma$ whose fiber $N_x\Sigma$ at each point $x \in
M$ is spanned by the vectors $(\mathbf{e}_3(x), \mathbf{e}_4(x))$ of any $3$-adapted frame at~$x$, together with
a~well-defined Riemannian metric
\begin{gather*}
I_{\rm normal} = \left(\omega^3_0\right)^2 + \left(\omega^4_0\right)^2
\end{gather*}
on $N\Sigma$.
\end{teo}

In order to obtain more information about the centroaf\/f\/ine metric, consider the structure
equations~\eqref{structure-eqns} for the semi-basic forms $\omega^1_0$, $\omega^2_0$ on~$M$.
Based on our adaptations, it is straightforward to compute that
\begin{gather*}
d\omega^1_0 = -\alpha \wedge \omega^2_0,
\qquad
d\omega^2_0 = \alpha \wedge \omega^1_0.
\end{gather*}
Therefore,~$\alpha$ is the Levi-Civita connection form associated to the centroaf\/f\/ine metric on~$\Sigma$, and the Gauss
curvature~$K$ of this metric is determined by the equation
\begin{gather*}
d\alpha = K  \omega^1_0 \wedge \omega^2_0.
\end{gather*}
The remaining structure equations~\eqref{structure-eqns} determine relations between the functions $h^i_{jk}$ on
$\mathcal{F}_3$ and their covariant derivatives with respect to $\omega^1_0$, $\omega^2_0$.
These relations may be viewed as analogs of the Gauss and Codazzi equations for Riemannian surfaces in Euclidean space.
In particular, the analog of the Gauss equation is
\begin{gather}
K = \tfrac{1}{2}\left(-h^3_{32} h^4_{31} - h^4_{41} h^3_{42} + h^3_{41} h^4_{42} - h^3_{32} h^4_{42} + h^4_{32}
h^3_{31} \right.
\nonumber
\\
\left.\phantom{K=}{}
+ h^4_{32} h^3_{42} - h^4_{41} h^3_{31} + h^3_{41} h^4_{31} + h^1_{31} - h^2_{32} + 2h^1_{42} \right) - \epsilon,
\label{space-like-Gauss-eqn}
\end{gather}
while the remainder of the relations are partial dif\/ferential equations involving the func\-tions~$h^i_{jk}$.
An analog of Bonnet's theorem (see~\cite{Griffiths74}) guarantees that, at least locally, any solution of this PDE
system gives rise to a~nondegenerate, space-like centroaf\/f\/ine surface, and that this surface is unique up to the action
of ${\rm GL}(5,\mathbb{R})$ on $\mathbb{R}^5 \setminus \{0 \}$.
In particular, the functions~$h^i_{jk}$ on~$\mathcal{F}_3$ form a~{\em complete} set of local invariants for such
surfaces.

\section{The time-like case}\label{time-like-sec}

Now, suppose that the plane~$P$ spanned by $(h^3, h^4)$ is time-like at every point of~$M$.
According to the group action~\eqref{plane-action}, we can f\/ind a~1-adapted frame f\/ield along~$\Sigma$ for which
$P=P_2$, as in equation~\eqref{normal-forms-planes}.
Furthermore, we can then use the action by~$B$ to f\/ind a~1-adapted frame f\/ield along~$\Sigma$ for which
\begin{gather}
\label{time-like-h3-h4}
h^3 =
\begin{bmatrix}
1 & 0
\\
0 & 0
\end{bmatrix},
\qquad
h^4 =
\begin{bmatrix}
0 & 0
\\
0 & 1
\end{bmatrix}.
\end{gather}
(This choice of $h^3$, $h^4$ represents a~null basis for $P_2$ with respect to the indef\/inite quadratic
form~\eqref{quadratic-form} on $P_2$.)

The next step is to determine the subgroup of $G_1$ that preserves the conditions~\eqref{time-like-h3-h4}.
To this end, f\/irst note that the plane $P_2$ spanned by the matrices~\eqref{time-like-h3-h4} consists precisely of the
diagonal matrices in ${\rm Sym}^2(\mathbb{R})$.
A~straightforward computation shows that, with $h^3$, $h^4$ as in~\eqref{time-like-h3-h4}, the of\/f-diagonal components of
$\tr(A^{\rm T} h^3 A)$ and $\tr(A^{\rm T} h^4 A)$ are equal to $a_{11}a_{12}$ and $a_{21}a_{22}$, respectively.
Therefore, the action~\eqref{plane-action} preserves $P_2$ if and only if
\begin{gather*}
a_{11}a_{12} = a_{21}a_{22} = 0.
\end{gather*}
Substituting this condition into equations~\eqref{1-adapted-group-action} and imposing the conditions $\tilde{h}^3 =
h^3$, $\tilde{h}^4 = h^4$ then shows that we must have either
\begin{gather}
A=
\begin{bmatrix}
a_{11} & 0
\\
0 & a_{22}
\end{bmatrix},
\qquad
B =
\begin{bmatrix}
a_{11}^2 & 0
\\
0 & a_{22}^2
\end{bmatrix},
\qquad
a_{11}, a_{22} \neq 0,
\label{time-like-AB-0}
\end{gather}
or
\begin{gather*}
A =
\begin{bmatrix}
0 & a_{12}
\\
a_{21} & 0
\end{bmatrix},
\qquad
B =
\begin{bmatrix}
0 & a_{12}^2
\\
a_{21}^2 & 0
\end{bmatrix},
\qquad
a_{12}, a_{21} \neq 0.
\end{gather*}
Since the latter transformation may be obtained from the former simply by interchanging $(\mathbf{e}_1, \mathbf{e}_2)$
and $(\mathbf{e}_3, \mathbf{e}_4)$, we will restrict our attention to transformations of the form~\eqref{time-like-AB-0},
where $A$, $B$ are diagonal matrices and $B = A^2$.

Next, consider the ef\/fect of the action~\eqref{1-adapted-group-action} on $h^0$.
With $A = I_2$ and $r_{03}$, $r_{04}$ chosen appropriately, we can add any linear combination of $h^3$, $h^4$ to $h^0$.
Thus we can f\/ind a~1-adapted frame f\/ield for which $h^0$ is a~multiple (nonzero by the nondegeneracy assumption) of the
matrix $
\begin{bmatrix}
0 & 1
\\
1 & 0
\end{bmatrix}
$.
Then under the action~\eqref{1-adapted-group-action} with $A$, $B$ as in~\eqref{time-like-AB-0}, we have
\begin{gather*}
\tilde{h}^0 = a_{11} a_{22} h^0 +
\begin{bmatrix}
-r_{03} & 0
\\
0 & -r_{04}
\end{bmatrix}.
\end{gather*}
Thus we can f\/ind a~1-adapted frame f\/ield along~$\Sigma$ satisfying the additional condition that
\begin{gather*}
h^0 =
\begin{bmatrix}
0 & 1
\\
1 & 0
\end{bmatrix},
\end{gather*}
and this condition is preserved by transformations of the form~\eqref{1-adapted-group-action} with $A$, $B$ as
in~\eqref{time-like-AB-0}, such that $a_{11} a_{22}=1$ and $r_{03} = r_{04} = 0$.
For simplicity, we will assume that $a_{11} > 0$; then we can set
\begin{gather*}
A=
\begin{bmatrix}
e^\lambda & 0
\\
0 & e^{-\lambda}
\end{bmatrix},
\qquad
B =
\begin{bmatrix}
e^{2\lambda} & 0
\\
0 & e^{-2\lambda}
\end{bmatrix},
\qquad
\lambda \in \mathbb{R}.
\end{gather*}

\begin{dfn}
Let $\Sigma = \bar{f}(M)$ be a~nondegenerate, time-like centroaf\/f\/ine surface in $\mathbb{R}^5 \setminus \{0\}$.
A~1-adapted frame f\/ield $f:M \to \mathcal{F}_1$ will be called {\em $2$-adapted} if it satisf\/ies the conditions
\begin{gather}
h^3 =
\begin{bmatrix}
1 & 0
\\
0 & 0
\end{bmatrix},
\qquad
h^4 =
\begin{bmatrix}
0 & 0
\\
0 & 1
\end{bmatrix},
\qquad
h^0 =
\begin{bmatrix}
0 & 1
\\
1 & 0
\end{bmatrix}
\label{2-adapted-time-like-h}
\end{gather}
at every point of~$M$.
\end{dfn}

Any two 2-adapted frames $(\mathbf{e}_0, \ldots, \mathbf{e}_4)$, $(\tilde{\mathbf{e}}_0, \ldots, \tilde{\mathbf{e}}_4)$
based at the same point $x \in M$ are related by a~transformation of the form
\begin{gather}
\begin{bmatrix}
\tilde{\mathbf{e}}_0 & \tilde{\mathbf{e}}_1 & \tilde{\mathbf{e}}_2 & \tilde{\mathbf{e}}_3 & \tilde{\mathbf{e}}_4
\end{bmatrix}
=
\begin{bmatrix}
\mathbf{e}_0 & \mathbf{e}_1 & \mathbf{e}_2 & \mathbf{e}_3 & \mathbf{e}_4
\end{bmatrix}
\begin{bmatrix}
1 & 0 & 0 & 0 & 0
\\
0 & e^{\lambda} & 0 & r_{13} & r_{14}
\\
0 & 0 & e^{-\lambda} & r_{23} & r_{24}
\\
0 & 0 & 0 & e^{2\lambda} & 0
\\
0 & 0 & 0 & 0 & e^{-2\lambda}
\end{bmatrix}.
\label{2-adapted-time-like-frame-transformation}
\end{gather}
We will denote the group of all matrices of the form in~\eqref{2-adapted-time-like-frame-transformation} by $G_2$; then
the 2-adapted frame f\/ields along~$\Sigma$ are the smooth sections of a~principal bundle $\mathcal{F}_2 \subset
\mathcal{F}_1$ over~$M$ with f\/iber group~$G_2$.

The equations~\eqref{2-adapted-time-like-h} are equivalent to the condition that the Maurer--Cartan forms associated to
a~2-adapted frame f\/ield satisfy the conditions
\begin{gather}
\omega^3_1=\omega^1_0,
\qquad
\omega^3_2=0,
\qquad
\omega^4_1=0,
\qquad
\omega^4_2=\omega^2_0,
\qquad
\omega^0_1=\omega^2_0,
\qquad
\omega^0_2=\omega^1_0.
\label{2-adapted-time-like-MC-relations}
\end{gather}
Dif\/ferentiating equations~\eqref{2-adapted-time-like-MC-relations} yields
\begin{alignat*}{3}
& \big(2\omega^1_1 - \omega^3_3\big) \wedge \omega^1_0 + \omega^1_2 \wedge \omega^2_0 = 0,
\qquad &&
\omega^1_2 \wedge \omega^1_0 - \omega^3_4 \wedge \omega^2_0 = 0,&
\\
& -\omega^4_3 \wedge \omega^1_0 + \omega^2_1 \wedge \omega^2_0 = 0,
\qquad &&
\omega^2_1 \wedge \omega^1_0 + \big(2\omega^2_2 - \omega^4_4\big) \wedge \omega^2_0 = 0,&
\\
& \big(2\omega^2_1 - \omega^0_3\big) \wedge \omega^1_0 + \big(\omega^1_1 + \omega^2_2\big) \wedge \omega^2_0 = 0,
\qquad &&
\big(\omega^1_1 + \omega^2_2\big) \wedge \omega^1_0 + \big(2\omega^1_2 - \omega^0_4\big) \wedge \omega^2_0 = 0.&
\end{alignat*}
Applying Cartan's lemma to these equations shows that there exists a~1-form~$\alpha$ and func\-tions~$h^i_{jk}$ on
$\mathcal{F}_2$ such that
\begin{alignat}{3}
&\omega^0_3=h^0_{31}\omega^1_0 + h^0_{32}\omega^2_0,
\qquad&&
\omega^0_4=h^0_{41}\omega^1_0 + h^0_{42}\omega^2_0,&
\nonumber
\\
&\omega^1_1=\alpha + h^1_{11}\omega^1_0 + h^2_{22}\omega^2_0,
\qquad&&
\omega^1_2=h^1_{21}\omega^1_0 + h^1_{22}\omega^2_0,&
\nonumber
\\
&\omega^2_1=h^2_{11}\omega^1_0 + h^2_{12}\omega^2_0,
\qquad&&
\omega^2_2=-\alpha + h^1_{11}\omega^1_0 + h^2_{22}\omega^2_0,&
\nonumber
\\
&\omega^3_3=2\alpha + h^3_{31}\omega^1_0 + h^3_{32}\omega^2_0,
\qquad&&
\omega^3_4=h^3_{41}\omega^1_0 + h^3_{42}\omega^2_0,&
\nonumber
\\
&\omega^4_3=h^4_{31}\omega^1_0 + h^4_{32}\omega^2_0,
\qquad&&
\omega^4_4=-2\alpha + h^4_{41}\omega^1_0 + h^4_{42}\omega^2_0.&
\label{2-adapted-time-like-MC-semi-basic-forms}
\end{alignat}
Moreover, the functions $h^i_{jk}$ satisfy the relations
\begin{alignat}{3}
& 2 h^2_{22} - h^1_{21} - h^3_{32} = 0,
\qquad &&
h^1_{22} + h^3_{41} = 0,&
\nonumber
\\
&h^2_{11} + h^4_{32} = 0,
\qquad &&
2 h^1_{11} - h^2_{12} - h^4_{41} = 0,&
\nonumber
\\
& 2 h^1_{11} - 2 h^2_{12} + h^0_{32} = 0,
\qquad &&
2 h^2_{22} - 2 h^1_{21} + h^0_{41} = 0.&
\label{2-adapted-time-like-h-relations}
\end{alignat}

If $f, \tilde{f}:M \to \mathcal{F}_2$ are two 2-adapted frame f\/ields related by a~transformation of the
form~\eqref{2-adapted-time-like-frame-transformation}, then we can once again use equation~\eqref{how-MC-forms-transform}
to determine how the corresponding functions $h^i_{jk}$, $\tilde{h}^i_{jk}$ are related.
As in the space-like case, the most straightforward to compute are those correspon\-ding to the forms $\tilde{\omega}^0_3$, $\tilde{\omega}^0_4$.
These forms appear as the coef\/f\/icients of $\tilde{\mathbf{e}}_0 = \mathbf{e}_0$ in the equations~\eqref{define-MC-forms}
for $d\tilde{\mathbf{e}}_3$, $d\tilde{\mathbf{e}}_4$.
By applying equations~\eqref{2-adapted-time-like-frame-transformation}
and~\eqref{2-adapted-time-like-MC-semi-basic-forms}, one can show that
\begin{gather*}
\begin{bmatrix}
\tilde{h}^0_{31}
\vspace{1mm}\\
\tilde{h}^0_{41}
\end{bmatrix}
=
\begin{bmatrix}
e^{2\lambda} & 0
\\
0 & e^{-2\lambda}
\end{bmatrix}
\begin{bmatrix}
h^0_{31}
\vspace{1mm}\\
h^0_{41}
\end{bmatrix}
+
\begin{bmatrix}
r_{23}
\\
r_{24}
\end{bmatrix},
\qquad
\begin{bmatrix}
\tilde{h}^0_{32}
\vspace{1mm}\\
\tilde{h}^0_{42}
\end{bmatrix}
=
\begin{bmatrix}
e^{2\lambda} & 0
\\
0 & e^{-2\lambda}
\end{bmatrix}
\begin{bmatrix}
h^0_{32}
\vspace{1mm}\\
h^0_{42}
\end{bmatrix}
+
\begin{bmatrix}
r_{13}
\\
r_{14}
\end{bmatrix}.
\end{gather*}
Thus we can f\/ind a~2-adapted frame f\/ield along~$\Sigma$ satisfying the conditions that
\begin{gather*}
h^0_{31} = h^0_{32} = h^0_{41} = h^0_{42} = 0,
\end{gather*}
and these conditions are preserved by transformations of the form~\eqref{2-adapted-time-like-frame-transformation} with
$r_{13} = r_{14} = r_{23} = r_{24} = 0$.

\begin{dfn}
Let $\Sigma = \bar{f}(M)$ be a~nondegenerate, time-like centroaf\/f\/ine surface in $\mathbb{R}^5 \setminus \{0\}$.
A~2-adapted frame f\/ield $f:M \to \mathcal{F}_2$ will be called {\em $3$-adapted} if it satisf\/ies the conditions
\begin{gather}
h^0_{31} = h^0_{32} = h^0_{41} = h^0_{42} = 0
\label{3-adapted-time-like-h}
\end{gather}
at every point of~$M$.
\end{dfn}

Any two 3-adapted frames $(\mathbf{e}_0, \ldots, \mathbf{e}_4)$, $(\tilde{\mathbf{e}}_0, \ldots, \tilde{\mathbf{e}}_4)$
based at the same point $x \in M$ are related by a~transformation of the form
\begin{gather}
\begin{bmatrix}
\tilde{\mathbf{e}}_0 & \tilde{\mathbf{e}}_1 & \tilde{\mathbf{e}}_2 & \tilde{\mathbf{e}}_3 & \tilde{\mathbf{e}}_4
\end{bmatrix}
=
\begin{bmatrix}
\mathbf{e}_0 & \mathbf{e}_1 & \mathbf{e}_2 & \mathbf{e}_3 & \mathbf{e}_4
\end{bmatrix}
\begin{bmatrix}
1 & 0 & 0 & 0 & 0
\\
0 & e^{\lambda} & 0 & 0 & 0
\\
0 & 0 & e^{-\lambda} & 0 & 0
\\
0 & 0 & 0 & e^{2\lambda} & 0
\\
0 & 0 & 0 & 0 & e^{-2\lambda}
\end{bmatrix}.
\label{3-adapted-time-like-frame-transformation}
\end{gather}
We will denote the group of all matrices of the form in~\eqref{3-adapted-time-like-frame-transformation} by $G_3$; note
that $G_3 \cong {\rm SO}^+(1,1)$.
Then the 3-adapted frame f\/ields are the smooth sections of a~principal bundle $\mathcal{F}_3 \subset \mathcal{F}_2$
over~$M$ with f\/iber group $G_3$.

The equations~\eqref{3-adapted-time-like-h} are equivalent to the condition that the Maurer--Cartan forms associated to
a~3-adapted frame f\/ield satisfy the conditions
\begin{gather}
\omega^0_3 = \omega^0_4 = 0.
\label{3-adapted-MC-time-like-relations}
\end{gather}
Dif\/ferentiating equations~\eqref{3-adapted-MC-time-like-relations} yields
\begin{gather*}
\omega^2_3 \wedge \omega^1_0 + \omega^1_3 \wedge \omega^2_0 = \omega^2_4 \wedge \omega^1_0 + \omega^1_4 \wedge
\omega^2_0 = 0,
\end{gather*}
and applying Cartan's lemma shows that there exist functions $h^i_{jk}$ on $\mathcal{F}_3$ such that
\begin{alignat*}{3}
& \omega^1_3=h^1_{31}\omega^1_0 + h^1_{32}\omega^2_0,
\qquad&&
\omega^2_3=h^2_{31}\omega^1_0 + h^1_{31}\omega^2_0,&
\\
& \omega^1_4=h^1_{41}\omega^1_0 + h^1_{42}\omega^2_0,
\qquad&&
\omega^2_4=h^2_{41}\omega^1_0 + h^1_{41}\omega^2_0.&
\end{alignat*}
Moreover, on $\mathcal{F}_3$, the last two relations in equations~\eqref{2-adapted-time-like-h-relations} simplify to
\begin{gather*}
h^1_{11} - h^2_{12} = h^2_{22} - h^1_{21} = 0.
\end{gather*}

At this point, we have canonically associated to any nondegenerate, time-like centroaf\/f\/ine surface in $\mathbb{R}^5
\setminus \{0\}$ a~frame bundle $\mathcal{F}_3$ over~$M$ with f\/iber group isomorphic to ${\rm SO}^+(1,1)$.
Thus we have the following theorem:

\begin{teo}
\label{time-like-thm}
Let $\bar{f}:M \to \mathbb{R}^5 \setminus\{0\}$ be a~centroaffine immersion whose image $\Sigma = \bar{f}(M)$ is
a~nondegenerate, time-like centroaffine surface.
Then the pullbacks of the Maurer--Cartan forms on ${\rm GL}(5,\mathbb{R})$ to the bundle $\mathcal{F}_3$ of $3$-adapted
frames on~$\Sigma$ determine a~well-defined Lorentzian metric
\begin{gather*}
I = 2   \omega^1_0  \omega^2_0
\end{gather*}
on~$\Sigma$, called the   centroaffine metric.
Moreoever, there is a~well-defined ``centroaffine normal bundle'' $N\Sigma$ whose fiber $N_x\Sigma$ at each point $x \in M$
is spanned by the vectors $(\mathbf{e}_3(x), \mathbf{e}_4(x))$ of any $3$-adapted frame at~$x$, together with a~well-defined Lorentzian metric
\begin{gather*}
I_{\rm normal} = 2\omega^3_0 \omega^4_0
\end{gather*}
on $N\Sigma$.
\end{teo}

In order to obtain more information about the centroaf\/f\/ine metric, consider the structure
equations~\eqref{structure-eqns} for the semi-basic forms $\omega^1_0$, $\omega^2_0$ on~$M$.
Based on our adaptations, it is straightforward to compute that
\begin{gather*}
d\omega^1_0 = -\alpha \wedge \omega^1_0,
\qquad
d\omega^2_0 = \alpha \wedge \omega^2_0.
\end{gather*}
Therefore,~$\alpha$ is the Levi-Civita connection form associated to the centroaf\/f\/ine metric on~$\Sigma$, and the Gauss
curvature~$K$ of this metric is determined by the equation
\begin{gather*}
d\alpha = K  \omega^1_0 \wedge \omega^2_0.
\end{gather*}
As in the space-like case, the remaining structure equations~\eqref{structure-eqns} determine relations between the
functions $h^i_{jk}$ on $\mathcal{F}_3$ and their covariant derivatives with respect to $\omega^1_0$, $\omega^2_0$.
The analog of the Gauss equation is
\begin{gather}
K = h^3_{41} h^4_{32} - h^3_{32} h^4_{41} + \tfrac{1}{2}\left(h^1_{32} + h^2_{41} \right) - 1,
\label{time-like-Gauss-eqn}
\end{gather}
while the remainder of the relations are partial dif\/ferential equations involving the functions $h^i_{jk}$.
As in the space-like case, any solution of this PDE system locally gives rise to a~nondegenerate, time-like centroaf\/f\/ine
surface; this surface is unique up to the action of ${\rm GL}(5,\mathbb{R})$ on $\mathbb{R}^5 \setminus \{0 \}$, and the
functions $h^i_{jk}$ on $\mathcal{F}_3$ form a~complete set of local invariants for such surfaces.

\section{Homogeneous examples}\label{homog-sec}

The goal of this section is to give a~complete classif\/ication (up to the ${\rm GL}(5,\mathbb{R})$-action on~$\mathbb{R}^5 \setminus\{0\}$) of the homogeneous examples of space-like and time-like nondegenerate centroaf\/f\/ine surfaces
in~$\mathbb{R}^5 \setminus\{0\}$.
First we must def\/ine precisely what we mean by the term ``homogeneous''. Because any such surface $\Sigma = \bar{f}(M)$
has a~well-def\/ined Riemannian or Lorentzian metric, any symmetry of~$\Sigma$ must preserve the centroaf\/f\/ine metric
on~$\Sigma$ and hence must in fact be an {\em isometry} of~$\Sigma$ with its centroaf\/f\/ine metric.
Furthermore, if the group of symmetries of~$\Sigma$ acts transitively, then the centroaf\/f\/ine metric must have constant
Gauss curvature~$K$.
As is well-known, the maximal isometry group of any Riemannian or Lorentzian surface has dimension less than or equal to
three, and in the maximal case the isometry group acts transitively on the orthonormal frame bundle.
Thus we will use the following def\/inition:

\begin{dfn}
\label{homogeneous-def}
Let $\mathcal{F}_3$ be the bundle of 3-adapted frames along a~nondegenerate, space-like or time-like centroaf\/f\/ine surface
$\Sigma = \bar{f}(M)$ in $\mathbb{R}^5 \setminus \{0 \}$.
A~dif\/feomorphism $\phi:\mathcal{F}_3 \to \mathcal{F}_3$ is called a~{\em symmetry} of~$\Sigma$ if $\phi^*\Omega =
\Omega$; i.e., if~$\phi$ preserves the Maurer--Cartan forms on $\mathcal{F}_3$.
A~nondegenerate space-like or time-like centroaf\/f\/ine surface $\Sigma = \bar{f}(M)$ will be called {\em homogeneous} if the
group of symmetries of~$\Sigma$ is a~3-dimensional Lie group that acts transitively on $\mathcal{F}_3$.
\end{dfn}

\begin{rmk}
It might also be of interest to consider the slightly less restrictive assumption that~$\Sigma$ has a~2-dimensional
group of symmetries that acts transitively on the base manifold~$M$, but we will not consider this scenario here.
\end{rmk}

Our procedure for classifying the homogeneous examples is as follows.
Observe that if~$\Sigma$ is homogeneous, then all the structure functions $h^i_{jk}$ must be constant on the bundle
$\mathcal{F}_3$ of 3-adapted frames on~$\Sigma$.
When this condition is imposed, the structure equations~\eqref{structure-eqns} become algebraic relations among the
constants $h^i_{jk}$.
Given any solution to these relations, the structure equations~\eqref{structure-eqns} imply that the corresponding
Maurer--Cartan form $\Omega = [\omega^i_j]$ takes values in a~3-dimensional Lie algebra $\mathfrak{g}$ which is realized
explicitly as a~Lie subalgebra of $\mathfrak{gl}(5,\mathbb{R})$.
Thus~$\Omega$ is also the Maurer--Cartan form of the connected Lie group $G \subset {\rm GL}(5,\mathbb{R})$ generated~by
exponentiating $\mathfrak{g}$, and this equivalence of the Maurer--Cartan forms on $\mathcal{F}_3$ with those on~$G$
implies that $\mathcal{F}_3$ is a~homogeneous space for~$G$; indeed,~$G$ must be precisely the symmetry group that was
assumed to act transitively on $\mathcal{F}_3$.

Now, choose any point $\mathbf{f}_0 = (\mathbf{e}_0, \mathbf{e}_1, \mathbf{e}_2, \mathbf{e}_3, \mathbf{e}_4) \in
\mathcal{F}_3$.
Recall that we can view $\mathbf{f}_0$ as an element of ${\rm GL}(5,\mathbb{R})$.
The centroaf\/f\/ine surface~$\Sigma$ is equivalent via the ${\rm GL}(5,\mathbb{R})$-action to the surface
$\widetilde{\Sigma} = \mathbf{f}_0^{-1} \cdot \Sigma$, and the bundle $\widetilde{\mathcal{F}}_3$ of 3-adapted frames
over $\widetilde{\Sigma}$ is given by $\widetilde{\mathcal{F}}_3 = \mathbf{f}_0^{-1} \cdot \mathcal{F}_3$.
So without loss of generality, we may assume that $\mathbf{f}_0$ is the identity matrix $I_5$.
With this assumption, the tangent space $T_{\mathbf{f}_0}\mathcal{F}_3$ is equal to the Lie algebra $\mathfrak{g}$, and
$\mathcal{F}_3$ must in fact be equal to~$G$.
Finally,~$\Sigma$ is given by the image of~$G$ under the projection~\eqref{define-pi}.

In order to carry out this procedure, we consider the space-like and time-like cases separately.

\subsection{Space-like homogeneous examples}

From the adaptations of Section~\ref{space-like-sec}, the matrix $\Omega = [\omega^i_j]$ of Maurer--Cartan forms on
$\mathcal{F}_3$ may be written~as
\begin{gather}
\Omega =
\begin{bmatrix}
0 & 0 & 0 & 0 & 0
\\
0 & 0 & 1 & 0 & 0
\\
0 & -1 & 0 & 0 & 0
\\
0 & 0 & 0 & 0 & 2
\\
0 & 0 & 0 & -2 & 0
\end{bmatrix}
\alpha
+
\begin{bmatrix}
0 & \epsilon & 0 & 0 & 0
\\
1 & \tfrac{1}{2}(h^3_{31} + h^3_{42}) - h^4_{32} + h^4_{41} & \tfrac{1}{2}(h^3_{32} - h^3_{41}) & h^1_{31} &
h^1_{41}
\vspace{1mm}\\
0 & \tfrac{1}{2}(h^3_{32} - h^3_{41}) & \tfrac{1}{2}(h^3_{31} + h^3_{42}) & h^1_{32} & h^1_{42}
\vspace{1mm}\\
0 & 1 & 0 & h^3_{31} & h^3_{41}
\vspace{1mm}\\
0 & 0 & 1 & h^4_{31} & h^4_{41}
\end{bmatrix}
\omega^1_0
\nonumber
\\
\phantom{\Omega =}{}
+
\begin{bmatrix}
0 & 0 & \epsilon & 0 & 0
\\
0 & \tfrac{1}{2}(h^3_{32} - h^3_{41}) & \tfrac{1}{2}(h^3_{31} + h^3_{42}) & h^1_{32} & h^1_{42}
\vspace{1mm}\\
1 & \tfrac{1}{2}(h^3_{31} + h^3_{42}) & \tfrac{1}{2}(h^3_{32} - h^3_{41}) + h^4_{31} + h^4_{42} & h^2_{32} &
h^2_{42}
\vspace{1mm}\\
0 & 0 & -1 & h^3_{32} & h^3_{42}
\vspace{1mm}\\
0 & 1 & 0 & h^4_{32} & h^4_{42}
\end{bmatrix}
\omega^2_0,
\label{space-like-reduced-MC-form}
\end{gather}
while the structure equations~\eqref{structure-eqns} may be written as
\begin{gather}
d\Omega = - \Omega \wedge \Omega.
\label{matrix-structure-eqns}
\end{gather}
Substituting~\eqref{space-like-reduced-MC-form} into equation~\eqref{matrix-structure-eqns} and imposing the condition
that all the functions $h^i_{jk}$ are constant leads to a~system of 25 algebraic equations (some linear, some quadratic)
for the 14 unknown constants $h^i_{jk}$.
A~somewhat tedious, but straightforward, computation shows that this system has precisely two solutions, one for
$\epsilon=1$ and one for $\epsilon=-1$.
These are described in the following two examples.

\begin{exmpl}
\label{space-like-epsilon-plus-one}
When $\epsilon=1$, the unique solution to equation~\eqref{matrix-structure-eqns} with all $h^i_{jk}$ constant is
\begin{gather}
\Omega =
\begin{bmatrix}
0 & 0 & 0 & 0 & 0
\\
0 & 0 & 1 & 0 & 0
\\
0 & -1 & 0 & 0 & 0
\\
0 & 0 & 0 & 0 & 2
\\
0 & 0 & 0 & -2 & 0
\end{bmatrix}
\alpha +
\begin{bmatrix}
0 & 1 & 0 & 0 & 0
\\
1 & 0 & 0 & \tfrac{1}{3} & 0
\\
0 & 0 & 0 & 0 & \tfrac{1}{3}
\\
0 & 1 & 0 & 0 & 0
\\
0 & 0 & 1 & 0 & 0
\end{bmatrix}
\omega^1_0 +
\begin{bmatrix}
0 & 0 & 1 & 0 & 0
\\
0 & 0 & 0 & 0 & \tfrac{1}{3}
\\
1 & 0 & 0 & -\tfrac{1}{3} & 0
\\
0 & 0 & -1 & 0 & 0
\\
0 & 1 & 0 & 0 & 0
\end{bmatrix}
\omega^2_0.
\label{space-like-epsilon-plus-one-MC-form}
\end{gather}
Furthermore, the Gauss equation~\eqref{space-like-Gauss-eqn} implies that the centroaf\/f\/ine metric has Gauss curvature $K
= -\tfrac{1}{3}$.

Denote the matrices in equation~\eqref{space-like-epsilon-plus-one-MC-form} by $M_0$, $M_1$, $M_2$, respectively, so that
\begin{gather*}
\Omega = M_0  \alpha + M_1  \omega^1_0 + M_2  \omega^2_0.
\end{gather*}
It is straightforward to compute that
\begin{gather*}
[M_0, M_1] = -M_2,
\qquad
[M_1, M_2] = \tfrac{1}{3} M_0,
\qquad
[M_2, M_0] = -M_1.
\end{gather*}
These bracket relations imply that the Lie algebra $\mathfrak{g} \subset \mathfrak{gl}(5,\mathbb{R})$ spanned by $(M_0,
M_1, M_2)$ is isomorphic to $\mathfrak{so}(1,2)$.
(They also suggest that a~more natural basis might be obtained by multiplying each of $M_1$, $M_2$ by $\sqrt{3}$.)
Furthermore, it is straightforward to check that $\mathfrak{g}$ acts irreducibly on $\mathbb{R}^5 \setminus \{0 \}$.
It is well-known (see, e.g.,~\cite[p.~149]{FultonHarris}) that $\mathfrak{so}(1,2)$ has a~unique irreducible
5-dimensional representation and that this representation arises from a~(unique) irreducible representation of ${\rm
SO}^+(1,2)$; it follows that the Lie group $G \subset {\rm GL}(5,\mathbb{R})$ corresponding to the Lie algebra
$\mathfrak{g}$ is isomorphic to ${\rm SO}^+(1,2)$.

The easiest way to compute a~local parametrization for~$G$~-- and hence for~$\Sigma$~-- is to compute the 1-parameter
subgroups generated by $M_0$, $M_1$, $M_2$ and take products of the resulting group elements.

{\bf Warning.} This must be done carefully in order to ensure that the resulting products cover the entire group~$G$,
which in turn guarantees that the resulting parametrization is surjective onto~$\Sigma$.
The subtlety of this issue can already be seen in the standard representation for $\mathfrak{so}(1,2)$: the basis
\begin{gather*}
\bar{M}_0 =
\begin{bmatrix}
0 & 0 & 0
\\
0 & 0 & 1
\\
0 & -1 & 0
\end{bmatrix},
\qquad
\bar{M}_1 =
\begin{bmatrix}
0 & 1 & 0
\\
1 & 0 & 0
\\
0 & 0 & 0
\end{bmatrix},
\qquad
\bar{M}_2 =
\begin{bmatrix}
0 & 0 & 1
\\
0 & 0 & 0
\\
1 & 0 & 0
\end{bmatrix}
\end{gather*}
has the same bracket relations as $(M_0, \sqrt{3}M_1, \sqrt{3}M_2)$, and exponentiating this basis yields the
1-parameter subgroups
\begin{gather*}
\bar{g}_0(t) =
\begin{bmatrix}
1 & 0 & 0
\\
0 & \cos(t) & \sin(t)
\\
0 & -\sin(t) & \cos(t)
\end{bmatrix},
\\
\bar{g}_1(u) =
\begin{bmatrix}
\cosh (u) & \sinh (u) & 0
\\
\sinh (u) & \cosh (u) & 0
\\
0 & 0 & 1
\end{bmatrix},
\qquad
\bar{g}_2(v) =
\begin{bmatrix}
\cosh (v) & 0 & \sinh (v)
\\
0 & 1 & 0
\\
\sinh (v) & 0 & \cosh (v)
\end{bmatrix}.
\end{gather*}
Now consider the following two maps $f_1, f_2:\mathbb{R}^3 \to {\rm SO}^+(1,2)$, which are obtained by multiplying the
elements $\bar{g}_0(t)$, $\bar{g}_1(u)$, $\bar{g}_2(v)$ in dif\/ferent orders:
\begin{gather*}
f_1(u,v,t) = \bar{g}_1(u) \bar{g}_2(v) \bar{g}_0(t)
\\
\hphantom{f_1(u,v,t)}{}
=
\begin{bmatrix}
\cosh(u) \cosh(v) & \begin{array}{@{}c@{}} \sinh(u) \cos(t) \vspace{-1mm}\\ - \cosh(u) \sinh(v) \sin(t)\end{array}   &
\begin{array}{@{}c@{}} \sinh(u) \sin(t) \vspace{-1mm}\\ + \cosh(u) \sinh(v) \cos(t)\end{array}
\vspace{1mm}\\
\sinh(u) \cosh(v) & \begin{array}{@{}c@{}}\cosh(u) \cos(t)\vspace{-1mm}\\- \sinh(u) \sinh(v) \sin(t)\end{array} &
\begin{array}{@{}c@{}}\cosh(u) \sin(t)\vspace{-1mm}\\ + \sinh(u) \sinh(v) \cos(t)\end{array}
\vspace{1mm}\\
\sinh(v) & -\cosh(v) \sin(t) & \cosh(v) \cos(t)
\end{bmatrix},
\\
f_2(u,v,t)=\bar{g}_1(u) \bar{g}_0(t) \bar{g}_2(v)
\\
\hphantom{f_2(u,v,t)}{}
=
\begin{bmatrix}
\begin{array}{@{}c@{}} \cosh(u) \cosh(v)\vspace{-1mm}\\+ \sinh(u) \sinh(v) \sin(t)\end{array} & \sinh(u) \cos(t) &
 \begin{array}{@{}c@{}} \cosh(u) \sinh(v)\vspace{-1mm}\\+ \sinh(u) \cosh(v) \sin(t)\end{array}
\vspace{1mm}\\
\begin{array}{@{}c@{}} \sinh(u) \cosh(v)\vspace{-1mm}\\+ \cosh(u) \sinh(v) \sin(t)\end{array}
 & \cosh(u) \cos(t) & \begin{array}{@{}c@{}} \sinh(u) \sinh(v)\vspace{-1mm}\\+ \cosh(u) \cosh(v) \sin(t)\end{array}
\vspace{1mm}\\
\sinh(v) \cos(t) & -\sin(t) & \cosh(v) \cos(t)
\end{bmatrix}.
\end{gather*}
It is not dif\/f\/icult to show that $f_1$ is surjective onto ${\rm SO}^+(1,2)$, whereas we can see from the middle column
that $f_2$ is not.
Thus we must perform this construction with care.

Now, the obvious correspondence
\begin{gather*}
\bar{M}_0 \leftrightarrow M_0,
\qquad
\bar{M}_1 \leftrightarrow \sqrt{3} M_1,
\qquad
\bar{M}_2 \leftrightarrow \sqrt{3} M_2,
\end{gather*}
def\/ines a~Lie algebra isomorphism between the standard representation of $\mathfrak{so}(1,2)$ and our Lie algebra
$\mathfrak{g} \subset \mathfrak{gl}(5,\mathbb{R})$.
Therefore, the surjectivity of the map $f_1$ above implies that the analogous map $f:\mathbb{R}^3 \to G$ will also be
surjective onto~$G$, and it follows that the map $\bar{f} = \pi \circ f$ will be surjective onto~$\Sigma$.
The map $\bar{f}$ will also turn out to be independent of~$t$, and when regarded as a~function of the two variables
$(u,v)$ it will def\/ine a~surjective parametrization $\bar{f}:\mathbb{R}^2 \to \Sigma$.

With these considerations in mind, def\/ine the 1-parameter subgroups
\begin{gather*}
g_0(t)= \exp(tM_0) =
\begin{bmatrix}
1 & 0 & 0 & 0 & 0
\\
0 & \cos(t) & \sin(t) & 0 & 0
\\
0 & -\sin(t) & \cos(t) & 0 & 0
\\
0 & 0 & 0 & \cos(2t) & \sin(2t)
\\
0 & 0 & 0 & -\sin(2t) & \cos(2t)
\end{bmatrix},
\\
g_1(u)= \exp(u \sqrt{3} M_1)
\\
\phantom{g_1(u)}
=
\begin{bmatrix}
\tfrac{1}{4}(3 \cosh(2u) + 1) & \tfrac{\sqrt{3}}{2}\sinh(2u) & 0 & \tfrac{1}{4} (\cosh(2u) - 1) & 0
\\
\tfrac{\sqrt{3}}{2}\sinh(2u) & \cosh(2u) & 0 & \tfrac{1}{2\sqrt{3}} \sinh(2u) & 0
\\
0 & 0 & \cosh(u) & 0 & \tfrac{1}{\sqrt{3}} \sinh(u)
\\
\tfrac{3}{4}(\cosh(2u)-1) & \tfrac{\sqrt{3}}{2}\sinh(2u) & 0 & \tfrac{1}{4}(\cosh(2u) + 3) & 0
\\
0 & 0 & \sqrt{3} \sinh(u) & 0 & \cosh(u)
\end{bmatrix},
\\
g_2(v)= \exp(v \sqrt{3} M_2)
\\
\phantom{g_2(v)}
=
\begin{bmatrix}
\tfrac{1}{4}(3 \cosh(2v) + 1) & 0 & \tfrac{\sqrt{3}}{2}\sinh(2v) & \tfrac{1}{4} (1 - \cosh(2v)) & 0
\\
0 & \cosh(v) & 0 & 0 & \tfrac{1}{\sqrt{3}}\sinh(v)
\\
\tfrac{\sqrt{3}}{2}\sinh(2v) & 0 & \cosh(2v) & -\tfrac{1}{2\sqrt{3}} \sinh(2v) & 0
\\
\tfrac{3}{4}(1 - \cosh(2v)) & 0 & -\tfrac{\sqrt{3}}{2}\sinh(2v) & \tfrac{1}{4}(\cosh(2v) + 3) & 0
\\
0 & \sqrt{3} \sinh(v) & 0 & 0 & \cosh(v)
\end{bmatrix}.
\end{gather*}
Then set
\begin{gather}
\bar{f}(u,v,t) = \pi\left(g_1(u) \cdot g_2(v) \cdot g_0(t) \right) =
\begin{bmatrix}
\tfrac{1}{2}\left(3\cosh^2(u)\cosh^2(v) - 1 \right)
\vspace{1mm}\\
\sqrt{3} \sinh(u) \cosh(u) \cosh^2(v)
\vspace{1mm}\\
\sqrt{3} \cosh(u) \sinh(v) \cosh(v)
\vspace{1mm}\\
\tfrac{3}{2}\left(\cosh^2(v) (\cosh^2(u) - 2) + 1 \right)
\vspace{1mm}\\
3 \sinh(u) \sinh(v) \cosh(v)
\end{bmatrix}.
\label{space-like-epsilon-plus-one-param}
\end{gather}
It follows from the discussion above that $\bar{f}$ is a~surjective map onto~$\Sigma$.
Moreover, it is straightforward to check that the tangent vectors $\bar{f}_u$, $\bar{f}_v$ are linearly independent for
all $(u,v) \in \mathbb{R}^2$; therefore $\bar{f}$ parametrizes a~smooth surface $\Sigma \subset \mathbb{R}^5 \setminus
\{0 \}$, as expected.

We can compute the centroaf\/f\/ine metric on~$\Sigma$ explicitly as follows.
Let $f:\mathbb{R}^3 \to G$ be the map corresponding to~\eqref{space-like-epsilon-plus-one-param}; i.e.,
\begin{gather*}
f(u,v,t) = g_1(u) \cdot g_2(v) \cdot g_0(t).
\end{gather*}
Then we have
\begin{gather*}
\Omega = M_0  \alpha + M_1  \omega^1_0 + M_2  \omega^2_0 = f^{-1}  df.
\end{gather*}
Comparing these two expressions for~$\Omega$ shows that
\begin{gather*}
\omega^1_0 = \sqrt{3}\left(\cosh (v) \cos(t)   du - \sin(t)  dv \right),
\qquad
\omega^2_0 = \sqrt{3}\left(\cosh (v) \sin(t)  du + \cos(t)  dv \right).
\end{gather*}
Therefore, the centroaf\/f\/ine metric on~$\Sigma$ is given~by
\begin{gather*}
I = \left(\omega^1_0\right)^2 + \left(\omega^2_0\right)^2 = 3\big(\cosh^2(v)  du^2 + dv^2\big).
\end{gather*}
As expected, this metric has constant Gauss curvature $K = -\frac{1}{3}$.
Topologically, $\Sigma$ is simply connected and isometric to the hyperbolic plane $\mathbb{H}$ of ``radius'' $\sqrt{3}$,
while $\mathcal{F}_3$ is isomorphic to the orthonormal frame bundle of $\mathbb{H}$, with~$\pi$ as the projection map.

In fact, we can describe~$\Sigma$ more intrinsically: if we denote the coordinates of a~point $\mathbf{x} \in
\mathbb{R}^5 \setminus \{0 \}$ as $(x_0, \ldots, x_4)$, then the coordinates of $\bar{f}(u,v)$ satisfy the quadratic
equations
\begin{gather}
x_1(x_0 - x_3 - 1) - x_2 x_4 = 0,
\qquad
x_2 (x_0 + x_3 - 1) - x_1 x_4 = 0,
\nonumber
\\
(4 x_0 - 1)^2 - 12 x_1^2 - 12 x_2^2 - 9 = 0.
\label{space-like-epsilon-plus-one-implicit-eqns}
\end{gather}
Therefore,~$\Sigma$ is contained in the (real) algebraic variety $X \subset \mathbb{R}^5 \setminus \{0 \} $ def\/ined~by
equations~\eqref{space-like-epsilon-plus-one-implicit-eqns}.
$\Sigma$ is not, however, equal to all of~$X$; for
instance,~$X$ contains an af\/f\/ine plane consisting of all points of the form $(1, 0, 0, x_3, x_4)$, and this plane
intersects~$\Sigma$ only when $x_3 = x_4 = 0$.
Projections of~$\Sigma$ to the $(x_1, x_2, x_0)$, $(x_1, x_2, x_3)$, and $(x_1, x_2, x_4)$ coordinate 3-planes are shown
in Fig.~\ref{space-like-epsilon-plus-one-figure}.
\begin{figure}[t] \centering
\includegraphics[width=2in]{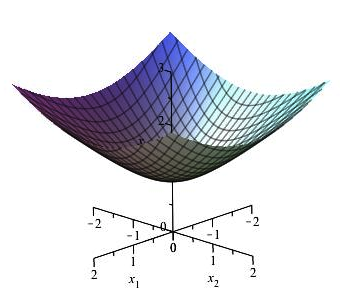} \includegraphics[width=2in]{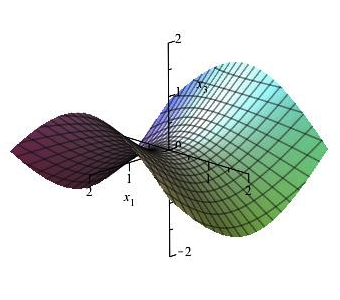} \includegraphics[width=2in]{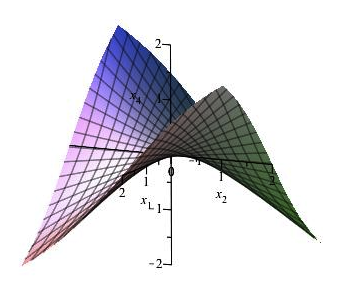}
\caption{Projections of surface from Example~\ref{space-like-epsilon-plus-one} to 3-D subspaces.}
\label{space-like-epsilon-plus-one-figure}
\end{figure}
\end{exmpl}

\begin{exmpl}
\label{space-like-epsilon-minus-one}
When $\epsilon=-1$, the unique solution to equation~\eqref{matrix-structure-eqns} with all $h^i_{jk}$ constant is
\begin{gather}
\Omega =
\begin{bmatrix}
0 & 0 & 0 & 0 & 0
\\
0 & 0 & 1 & 0 & 0
\\
0 & -1 & 0 & 0 & 0
\\
0 & 0 & 0 & 0 & 2
\\
0 & 0 & 0 & -2 & 0
\end{bmatrix}
\alpha +
\begin{bmatrix}
0 & -1 & 0 & 0 & 0
\\
1 & 0 & 0 & -\tfrac{1}{3} & 0
\\
0 & 0 & 0 & 0 & -\tfrac{1}{3}
\\
0 & 1 & 0 & 0 & 0
\\
0 & 0 & 1 & 0 & 0
\end{bmatrix}
\omega^1_0 +
\begin{bmatrix}
0 & 0 & -1 & 0 & 0
\\
0 & 0 & 0 & 0 & -\tfrac{1}{3}
\\
1 & 0 & 0 & \tfrac{1}{3} & 0
\\
0 & 0 & -1 & 0 & 0
\\
0 & 1 & 0 & 0 & 0
\end{bmatrix}
\omega^2_0.
\label{space-like-epsilon-minus-one-MC-form}
\end{gather}
Furthermore, the Gauss equation~\eqref{space-like-Gauss-eqn} implies that the centroaf\/f\/ine metric has Gauss curvature $K
= \tfrac{1}{3}$.

As in the previous example, denote the matrices in equation~\eqref{space-like-epsilon-minus-one-MC-form} by $M_0$, $M_1$, $M_2$, respectively.
Then we have
\begin{gather*}
[M_0, M_1] = -M_2,
\qquad
[M_1, M_2] = -\tfrac{1}{3} M_0,
\qquad
[M_2, M_0] = -M_1.
\end{gather*}
These bracket relations imply that the Lie algebra $\mathfrak{g} \subset \mathfrak{gl}(5,\mathbb{R})$ spanned by $(M_0,
M_1, M_2)$ is isomorphic to $\mathfrak{so}(3,\mathbb{R})$.
Furthermore, it is straightforward to check that $\mathfrak{g}$ acts irreducibly on~$\mathbb{R}^5 \setminus \{0 \}$.
Similarly to the previous example, $\mathfrak{so}(3,\mathbb{R})$ has a~unique irreducible 5-dimensional representation,
and this representation arises from a~(unique) irreducible representation of~${\rm SO}(3,\mathbb{R})$; it follows that
the Lie group $G \subset {\rm GL}(5,\mathbb{R})$ corresponding to the Lie algebra $\mathfrak{g}$ is isomorphic to~${\rm
SO}(3,\mathbb{R})$.

We will compute a~local parametrization for~$\Sigma$ as in the previous example: compute the 1-parameter subgroups
of~$G$ generated by $M_0$, $M_1$, $M_2$ and take products of the resulting group elements.
Ensuring the surjectivity of the resulting parametrization is easier than in the previous example.
First, observe that a~basis $(\bar{M}_0, \bar{M}_1, \bar{M}_2)$ for the standard representation of
$\mathfrak{so}(3,\mathbb{R})$ with the same bracket relations as $(M_0, \sqrt{3}M_1, \sqrt{3}M_2)$ is given~by
\begin{gather*}
\bar{M}_0 =
\begin{bmatrix}
0 & 0 & 0
\\
0 & 0 & 1
\\
0 & -1 & 0
\end{bmatrix},
\qquad
\bar{M}_1 =
\begin{bmatrix}
0 & -1 & 0
\\
1 & 0 & 0
\\
0 & 0 & 0
\end{bmatrix},
\qquad
\bar{M}_2 =
\begin{bmatrix}
0 & 0 & -1
\\
0 & 0 & 0
\\
1 & 0 & 0
\end{bmatrix}.
\end{gather*}
Exponentiating this basis yields the 1-parameter subgroups
\begin{gather*}
\bar{g}_0(t) =
\begin{bmatrix}
1 & 0 & 0
\\
0 & \cos(t) & \sin(t)
\\
0 & -\sin(t) & \cos(t)
\end{bmatrix},
\\
\bar{g}_1(u) =
\begin{bmatrix}
\cos (u) & -\sin (u) & 0
\\
\sin (u) & \cos (u) & 0
\\
0 & 0 & 1
\end{bmatrix},
\qquad
\bar{g}_2(v) =
\begin{bmatrix}
\cos (v) & 0 & -\sin (v)
\\
0 & 1 & 0
\\
\sin (v) & 0 & \cos (v)
\end{bmatrix}.
\end{gather*}
Then the map $f: \mathbb{R}^3 \to {\rm SO}(3,\mathbb{R})$ def\/ined~by
\begin{gather*}
f(u,v,t) = \bar{g}_1(u) \bar{g}_2(v) \bar{g}_0(t)
\\
\hphantom{f(u,v,t)}{} =
\begin{bmatrix}
\cos(u) \cos(v) & \begin{array}{@{}c@{}} -\sin(u) \cos(t)\vspace{-1mm}\\+ \cos(u) \sin(v) \sin(t)\end{array} & \begin{array}{@{}c@{}} -\sin(u) \sin(t)\vspace{-1mm}\\ - \cos(u) \sin(v) \cos(t)\end{array}
\vspace{1mm}\\
\sin(u) \cos(v) & \begin{array}{@{}c@{}} \cos(u) \cos(t)\vspace{-1mm}\\ + \sin(u) \sin(v) \sin(t)\end{array} & \begin{array}{@{}c@{}} \cos(u) \sin(t)\vspace{-1mm}\\- \sin(u) \sin(v) \cos(t)\end{array}
\vspace{1mm}\\
\sin(v) & -\cos(v) \sin(t) & \cos(v) \cos(t)
\end{bmatrix}
\end{gather*}
is easily seen to be surjective onto ${\rm SO}(3,\mathbb{R})$.
Thus the analogous map $f:\mathbb{R}^3 \to G$ will be surjective onto~$G$, and the map $\bar{f} = \pi \circ f$ will be
surjective onto~$\Sigma$.

So, def\/ine the 1-parameter subgroups
\begin{gather*}
g_0(t)=\exp(tM_0) =
\begin{bmatrix}
1 & 0 & 0 & 0 & 0
\\
0 & \cos(t) & \sin(t) & 0 & 0
\\
0 & -\sin(t) & \cos(t) & 0 & 0
\\
0 & 0 & 0 & \cos(2t) & \sin(2t)
\\
0 & 0 & 0 & -\sin(2t) & \cos(2t)
\end{bmatrix},
\\
g_1(u)=\exp(u \sqrt{3} M_1)
\\
\phantom{g_1(u)}
=
\begin{bmatrix}
\tfrac{1}{4}(3 \cos(2u) + 1) & -\tfrac{\sqrt{3}}{2}\sin(2u) & 0 & \tfrac{1}{4} (1 - \cos(2u)) & 0
\\
\tfrac{\sqrt{3}}{2}\sin(2u) & \cos(2u) & 0 & -\tfrac{1}{2\sqrt{3}} \sin(2u) & 0
\\
0 & 0 & \cos(u) & 0 & -\tfrac{1}{\sqrt{3}} \sin(u)
\\
\tfrac{3}{4}(1 - \cos(2u)) & \tfrac{\sqrt{3}}{2}\sin(2u) & 0 & \tfrac{1}{4}(\cos(2u) + 3) & 0
\\
0 & 0 & \sqrt{3} \sin(u) & 0 & \cos(u)
\end{bmatrix},
\\
g_2(v)=\exp(v \sqrt{3} M_2)
\\
\phantom{g_2(v)}
=
\begin{bmatrix}
\tfrac{1}{4}(3 \cos(2v) + 1) & 0 & -\tfrac{\sqrt{3}}{2}\sin(2v) & \tfrac{1}{4} (\cos(2v) - 1) & 0
\\
0 & \cos(v) & 0 & 0 & -\tfrac{1}{\sqrt{3}}\sin(v)
\\
\tfrac{\sqrt{3}}{2}\sin(2v) & 0 & \cos(2v) & \tfrac{1}{2\sqrt{3}} \sin(2v) & 0
\\
\tfrac{3}{4}(\cos(2v) - 1) & 0 & -\tfrac{\sqrt{3}}{2}\sin(2v) & \tfrac{1}{4}(\cos(2v) + 3) & 0
\\
0 & \sqrt{3} \sin(v) & 0 & 0 & \cos(v)
\end{bmatrix}.
\end{gather*}
Then set
\begin{gather}
\bar{f}(u,v,t) = \pi\left(g_1(u) \cdot g_2(v) \cdot g_0(t) \right) =
\begin{bmatrix}
\tfrac{1}{2}\left(3\cos^2(u)\cos^2(v) - 1 \right)
\vspace{1mm}\\
\sqrt{3} \sin(u) \cos(u) \cos^2(v)
\vspace{1mm}\\
\sqrt{3} \cos(u) \sin(v) \cos(v)
\vspace{1mm}\\
\tfrac{3}{2}\left(\cos^2(v) (2 - \cos^2(u)) - 1 \right)
\vspace{1mm}\\
3 \sin(u) \sin(v) \cos(v)
\end{bmatrix}.
\label{space-like-epsilon-minus-one-param}
\end{gather}
It follows from the discussion above that $\bar{f}$ is a~surjective map onto~$\Sigma$, and we see that $\bar{f}$ is also
independent of~$t$.
Unlike in the previous example, the tangent vectors $\bar{f}_u$, $\bar{f}_v$ are not linearly independent for all $(u,v)
\in \mathbb{R}^2$; indeed, $\bar{f}_u=0$ whenever~$v$ is an odd multiple of $\frac{\pi}{2}$.
Nevertheless, the restriction of $\bar{f}$ to some neighborhood of the point $(u,v) = (0,0)$ is a~smooth embedding, and
then homogeneity implies that~$\Sigma$ is smooth everywhere.

We can compute the centroaf\/f\/ine metric on~$\Sigma$ as in the previous example.
Let $f:\mathbb{R}^3 \to G$ be the map corresponding to~\eqref{space-like-epsilon-minus-one-param}; i.e.,
\begin{gather*}
f(u,v,t) = g_1(u) \cdot g_2(v) \cdot g_0(t).
\end{gather*}
Then we have
\begin{gather*}
\Omega = M_0  \alpha + M_1  \omega^1_0 + M_2  \omega^2_0 = f^{-1}  d f.
\end{gather*}
Comparing these two expressions for~$\Omega$ shows that, whenever $\cos(v) \neq 0$,
\begin{gather*}
\omega^1_0 = \sqrt{3}\big(\cos (v) \cos(t)  du - \sin(t) dv \big),
\qquad
\omega^2_0 = \sqrt{3}\big(\cos (v) \sin(t)   du + \cos(t)  dv \big).
\end{gather*}
Therefore, the centroaf\/f\/ine metric on~$\Sigma$ is given~by
\begin{gather*}
I = \big(\omega^1_0\big)^2 + \big(\omega^2_0\big)^2 = 3\big(\cos^2(v)  du^2 + dv^2\big).
\end{gather*}
As expected, this metric has constant Gauss curvature $K = \frac{1}{3}$.
Topologically, $\Sigma$ is simply connected and isometric to the sphere~$S$ of radius $\sqrt{3}$, while $\mathcal{F}_3$
is isomorphic to the orthonormal frame bundle of~$S$, with~$\pi$ as the projection map.

As in the previous example, we can show that~$\Sigma$ is contained in the intersection of three quadric hypersurfaces in
$\mathbb{R}^5 \setminus \{0 \}$.
The coordinates of $\bar{f}(u,v)$ satisfy the quadratic equations
\begin{gather}
x_1(x_0 + x_3 - 1) + x_2 x_4 = 0,
\qquad
x_2 (x_0 - x_3 - 1) + x_1 x_4 = 0,
\nonumber
\\
(4 x_0 - 1)^2 + 12 x_1^2 + 12 x_2^2 - 9 = 0.
\label{space-like-epsilon-minus-one-implicit-eqns}
\end{gather}
Therefore, $\Sigma$ is contained in the (real) algebraic variety $X \subset \mathbb{R}^5 \setminus \{0 \} $ def\/ined~by
equations~\eqref{space-like-epsilon-minus-one-implicit-eqns}.
Once again, $\Sigma$ is not equal to all of~$X$: the variety~$X$ contains af\/f\/ine planes consisting of all points of the
form $(1, 0, 0, x_3, x_4)$ or $(-\tfrac{1}{2}, 0, 0, x_3, x_4)$; the former intersects~$\Sigma$ only when $x_3 = x_4 =
0$, and the latter intersects~$\Sigma$ only when $x_3^2 + x_4^2 = \tfrac{9}{4}$.
Projections of~$\Sigma$ to the $(x_1, x_2, x_0)$, $(x_1, x_2, x_3)$, and $(x_1, x_2, x_4)$ coordinate 3-planes are shown
in Fig.~\ref{space-like-epsilon-minus-one-figure}.
\begin{figure}[h]\centering
\includegraphics[width=2in]{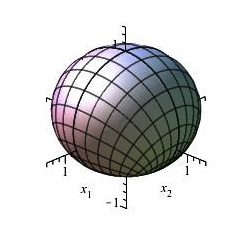} \includegraphics[width=2in]{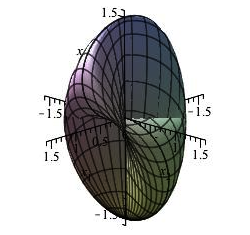} \includegraphics[width=2in]{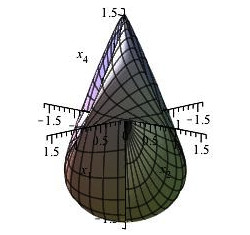}
\caption{Projections of surface from Example~\ref{space-like-epsilon-minus-one} to 3-D subspaces.}
\label{space-like-epsilon-minus-one-figure}
\end{figure}
\end{exmpl}

\subsection{Time-like homogeneous examples}

From the adaptations of Section~\ref{time-like-sec}, the matrix $\Omega = [\omega^i_j]$ of Maurer--Cartan forms on
$\mathcal{F}_3$ may be written~as
\begin{gather}
\Omega =
\begin{bmatrix}
0 & 0 & 0 & 0 & 0
\\
0 & 1 & 0 & 0 & 0
\\
0 & 0 & -1 & 0 & 0
\\
0 & 0 & 0 & 2 & 0
\\
0 & 0 & 0 & 0 & -2
\end{bmatrix}
\alpha +
\begin{bmatrix}
0 & 0 & 1 & 0 & 0
\\
1 & h^4_{41} & h^3_{32} & h^1_{31} & h^1_{41}
\vspace{1mm}\\
0 & -h^4_{32} & h^4_{41} & h^2_{31} & h^2_{41}
\vspace{1mm}\\
0 & 1 & 0 & h^3_{31} & h^3_{41}
\vspace{1mm}\\
0 & 0 & 0 & h^4_{31} & h^4_{41}
\end{bmatrix}
\omega^1_0 \nonumber\\
\hphantom{\Omega =}{}
+
\begin{bmatrix}
0 & 1 & 0 & 0 & 0
\\
0 & h^3_{32} & -h^3_{41} & h^1_{32} & h^1_{42}
\vspace{1mm}\\
1 & h^4_{41} & h^3_{32} & h^1_{31} & h^1_{41}
\vspace{1mm}\\
0 & 0 & 0 & h^3_{32} & h^3_{42}
\vspace{1mm}\\
0 & 0 & 1 & h^4_{32} & h^4_{42}
\end{bmatrix}
\omega^2_0,
\label{time-like-reduced-MC-form}
\end{gather}
while the structure equations~\eqref{structure-eqns} may be written as
\begin{gather}
d\Omega = - \Omega \wedge \Omega.
\label{time-like-matrix-structure-eqns}
\end{gather}
Substituting~\eqref{time-like-reduced-MC-form} into equation~\eqref{time-like-matrix-structure-eqns} and imposing the
condition that all the functions $h^i_{jk}$ are constant leads to a~system of 23 algebraic equations for the 14 unknown
constants $h^i_{jk}$.
This system has precisely one solution, which is described in the following example.

\begin{exmpl}
\label{time-like-example}
The unique solution to equation~\eqref{time-like-matrix-structure-eqns} with all $h^i_{jk}$ constant is
\begin{gather}
\Omega =
\begin{bmatrix}
0 & 0 & 0 & 0 & 0
\\
0 & 1 & 0 & 0 & 0
\\
0 & 0 & -1 & 0 & 0
\\
0 & 0 & 0 & 2 & 0
\\
0 & 0 & 0 & 0 & -2
\end{bmatrix}
\alpha +
\begin{bmatrix}
0 & 0 & 1 & 0 & 0
\\
1 & 0 & 0 & 0 & 0
\\
0 & 0 & 0 & 0 & \tfrac{2}{3}
\\
0 & 1 & 0 & 0 & 0
\\
0 & 0 & 0 & 0 & 0
\end{bmatrix}
\omega^1_0 +
\begin{bmatrix}
0 & 1 & 0 & 0 & 0
\\
0 & 0 & 0 & \tfrac{2}{3} & 0
\\
1 & 0 & 0 & 0 & 0
\\
0 & 0 & 0 & 0 & 0
\\
0 & 0 & 1 & 0 & 0
\end{bmatrix}
\omega^2_0.
\label{time-like-MC-form}
\end{gather}
Furthermore, the Gauss equation~\eqref{time-like-Gauss-eqn} implies that the centroaf\/f\/ine metric has Gauss curvature $K =
-\tfrac{1}{3}$.

As in the space-like examples, denote the matrices in equation~\eqref{time-like-MC-form} by $M_0$, $M_1$, $M_2$, respectively,
so that
\begin{gather*}
\Omega = M_0 \alpha + M_1  \omega^1_0 + M_2  \omega^2_0.
\end{gather*}
Then we have
\begin{gather*}
[M_0, M_1] = M_1,
\qquad
[M_1, M_2] = \tfrac{1}{3} M_0,
\qquad
[M_2, M_0] = M_2.
\end{gather*}
These bracket relations imply that the Lie algebra $\mathfrak{g} \subset \mathfrak{gl}(5,\mathbb{R})$ spanned by $(M_0,
M_1, M_2)$ is isomorphic to $\mathfrak{so}(2,1)$.
Furthermore, it is straightforward to check that $\mathfrak{g}$ acts irreducibly on~$\mathbb{R}^5 \setminus \{0 \}$.
Similarly to the previous examples, $\mathfrak{so}(2,1)$ has a~unique irreducible 5-dimensional representation, and this
representation arises from a~(unique) irreducible representation of ${\rm SO}^+(2,1)$; it follows that the Lie group $G
\subset {\rm GL}(5,\mathbb{R})$ corresponding to the Lie algebra $\mathfrak{g}$ is isomorphic to~${\rm SO}^+(2,1)$.

We will compute a~local parametrization for~$\Sigma$ more or less as in the space-like examples, by computing 1-parameter
subgroups of~$G$ and taking products of the resulting group elements.
Unfortunately, the basis $(M_0, M_1, M_2)$ is not well-suited to generating a~surjective parametrization, so f\/irst we
need to modify it slightly.
To this end, observe that a~basis $(\bar{M}_0, \bar{M}_1, \bar{M}_2)$ for the standard representation of
$\mathfrak{so}(2,1)$ with the same bracket relations as $(M_0, \sqrt{3}M_1, \sqrt{3}M_2)$ is given~by
\begin{gather*}
\bar{M}_0 =
\begin{bmatrix}
0 & 0 & 0
\\
0 & 1 & 0
\\
0 & 0 & -1
\end{bmatrix},
\qquad
\bar{M}_1 =
\begin{bmatrix}
0 & 0 & 1
\\
1 & 0 & 0
\\
0 & 0 & 0
\end{bmatrix},
\qquad
\bar{M}_2 =
\begin{bmatrix}
0 & 1 & 0
\\
0 & 0 & 0
\\
1 & 0 & 0
\end{bmatrix}.
\end{gather*}
Now consider the modif\/ied basis
\begin{gather*}
\bar{M}_0' = \bar{M}_0,
\qquad
\bar{M}_1' = \frac{1}{\sqrt{2}}(\bar{M}_1 - \bar{M}_2),
\qquad
\bar{M}_2' = \frac{1}{\sqrt{2}}(\bar{M}_1 + \bar{M}_2).
\end{gather*}
Exponentiating this modif\/ied basis yields the 1-parameter subgroups
\begin{gather*}
\bar{g}_0(t) =
\begin{bmatrix}
1 & 0 & 0
\\
0 & e^t & 0
\\
0 & 0 & e^{-t}
\end{bmatrix},
\qquad
\bar{g}_1(u) =
\begin{bmatrix}
\cos(u) & -\frac{1}{\sqrt{2}}\sin(u) & \frac{1}{\sqrt{2}}\sin(u)
\\
\frac{1}{\sqrt{2}}\sin(u) & \frac{1}{2}(1 + \cos(u)) & \frac{1}{2}(1 - \cos(u))
\\
-\frac{1}{\sqrt{2}}\sin(u) & \frac{1}{2}(1 - \cos(u)) & \frac{1}{2}(1 + \cos(u))
\end{bmatrix},
\\
\bar{g}_2(v) =
\begin{bmatrix}
\cosh(v) & \frac{1}{\sqrt{2}}\sinh(v) & \frac{1}{\sqrt{2}}\sinh(v)
\\
\frac{1}{\sqrt{2}}\sinh(v) & \frac{1}{2}(\cosh(v) + 1) & \frac{1}{2}(\cosh(v) - 1)
\\
\frac{1}{\sqrt{2}}\sinh(v) & \frac{1}{2}(\cosh(v) - 1) & \frac{1}{2}(\cosh(v) + 1)
\end{bmatrix}.
\end{gather*}
Then the map $f: \mathbb{R}^3 \to {\rm SO}^+(2,1)$ def\/ined~by
\begin{gather*}
f(u,v,t) = \bar{g}_1(u) \bar{g}_2(v) \bar{g}_0(t)\\
\hphantom{f(u,v,t)}{} = \begin{bmatrix}
\cos(u) \cosh(v) & \begin{array}{@{}c@{}}\frac{1}{\sqrt{2}} e^t(\cos(u) \sinh(v)\vspace{-1mm}\\ -\sin(u))\end{array} &
\begin{array}{@{}c@{}}\frac{1}{\sqrt{2}} e^{-t}(\cos(u) \sinh(v)\vspace{-1mm}\\+\sin(u))\end{array}
\vspace{1mm}\\
\begin{array}{@{}c@{}}\frac{1}{\sqrt{2}}(\sin(u)\cosh(v)\vspace{-1mm}\\+\sinh(v))\end{array} &
\begin{array}{@{}c@{}}\frac{1}{2}e^t(\sin(u)\sinh(v)\vspace{-1mm}\\+\cosh(v)+\cos(u))\end{array} &
\begin{array}{@{}c@{}}\frac{1}{2}e^{-t}(\sin(u) \sinh(v)\vspace{-1mm}\\ +\cosh(v)-\cos(u))\end{array}
\vspace{1mm}\\
\begin{array}{@{}c@{}}-\frac{1}{\sqrt{2}}(\sin(u)\cosh(v)\vspace{-1mm}\\ -\sinh(v))\end{array} &
\begin{array}{@{}c@{}}-\frac{1}{2}e^t(\sin(u)\sinh(v)\vspace{-1mm}\\ -\cosh(v)+\cos(u))\end{array} &
\begin{array}{@{}c@{}}-\frac{1}{2}e^{-t}(\sin(u)\sinh(v)\vspace{-1mm}\\ -\cosh(v)-\cos(u))\end{array}
\end{bmatrix}
\end{gather*}
is surjective onto ${\rm SO}^+(2,1)$.
Thus the analogous map $f:\mathbb{R}^3 \to G$ will be surjective onto~$G$, and the map $\bar{f} = \pi \circ f$ will be
surjective onto~$\Sigma$.

So, def\/ine the 1-parameter subgroups
\begin{gather*}
g_0(t) = \exp(tM_0),
\qquad\! \!
g_1(u) = \exp\Big(u \sqrt{\tfrac{3}{2}} (M_1 - M_2)\Big),
\qquad\! \!
g_2(v) = \exp\Big(v \sqrt{\tfrac{3}{2}} (M_1 + M_2)\Big).
\end{gather*}
(The explicit expressions for these group elements are each too large to f\/it on one line and are not particularly
enlightening.) Then set
\begin{gather}
\bar{f}(u,v,t)=\pi\left(g_1(u) \cdot g_2(v) \cdot g_0(t) \right)
\nonumber
\\
\phantom{\bar{f}(u,v,t)}
  =
\begin{bmatrix}
\frac{1}{4}[3 \cos^2(u)(\cosh(2v) + 1) - 2]
\vspace{1mm}\\
\frac{\sqrt{6}}{4} \cos(u) [\sin(u)(\cosh(2v) + 1) + \sinh(2v)]
\vspace{1mm}\\
-\frac{\sqrt{6}}{4} \cos(u) [\sin(u)(\cosh(2v) + 1) - \sinh(2v)]
\vspace{1mm}\\
-\frac{3}{8} [\cos^2(u)(\cosh(2v) + 1) - 2(\cosh(2v) + \sin(u) \sinh(2v))]
\vspace{1mm}\\
-\frac{3}{8} [\cos^2(u)(\cosh(2v) + 1) - 2(\cosh(2v) - \sin(u) \sinh(2v))]
\end{bmatrix}.
\label{time-like-param}
\end{gather}
It follows from the discussion above that $\bar{f}$ is a~surjective map onto~$\Sigma$, and we see that $\bar{f}$ is also
independent of~$t$.
Moreover, it is straightforward to check that the tangent vectors $\bar{f}_u$, $\bar{f}_v$ are linearly independent for
all $(u,v) \in \mathbb{R}^2$; therefore $\bar{f}$ parametrizes a~smooth surface $\Sigma \subset \mathbb{R}^5 \setminus
\{0 \}$, as expected.

We can compute the centroaf\/f\/ine metric on~$\Sigma$ as in the previous examples.
Let $f:\mathbb{R}^3 \to G$ be the map corresponding to~\eqref{time-like-param}; i.e.,
\begin{gather*}
f(u,v,t) = g_1(u) \cdot g_2(v) \cdot g_0(t).
\end{gather*}
Then we have
\begin{gather*}
\Omega = M_0  \alpha + M_1  \omega^1_0 + M_2  \omega^2_0 = f^{-1}  df.
\end{gather*}
Comparing these two expressions for~$\Omega$ shows that
\begin{gather*}
\omega^1_0 = \sqrt{\frac{3}{2}} e^{-t} \left(\cosh (v)  du + dv \right),
\qquad
\omega^2_0 = \sqrt{\frac{3}{2}} e^{t} \left(-\cosh (v)  du + dv \right).
\end{gather*}
Therefore, the centroaf\/f\/ine metric on~$\Sigma$ is given~by
\begin{gather*}
I = 2 \omega^1_0 \omega^2_0 = 3\big({-}\cosh^2(v)  du^2 + dv^2\big).
\end{gather*}
As expected, this metric has constant Gauss curvature $K = \frac{1}{3}$.

In this case, $\mathcal{F}_3$ is isomorphic to the orthonormal frame bundle of the time-like surface $S^{2,1}$ consisting
of all space-like vectors of length $\sqrt{3}$ in $\mathbb{R}^{2,1}$.
This surface is a~hyperboloid of one sheet, and so we might expect that~$\Sigma$ would be dif\/feomorphic to this
hyperboloid.
However, if we regard the domain of the parametrization~\eqref{time-like-param} as $S^1 \times \mathbb{R}$, we see that
the map~\eqref{time-like-param} is invariant under the transformation
\begin{gather}
(u,v) \to (u+\pi, -v).
\label{involution}
\end{gather}
Therefore,~$\Sigma$ is dif\/feomorphic to a~M\"obius band, which is precisely the quotient of the hyperboloid by this map.
(We note, however, that the map $f: S^1 \times \mathbb{R} \to G$ {\em is} one-to-one, as the frame vectors
$\mathbf{e}_1(u,v), \ldots, \mathbf{e}_4(u,v)$ are not preserved by the transformation~\eqref{involution}.)

\begin{figure}[t]\centering
\includegraphics[width=2in]{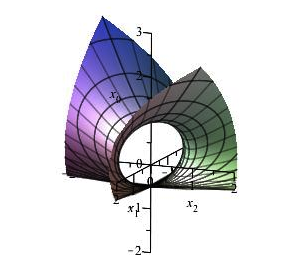} \includegraphics[width=2in]{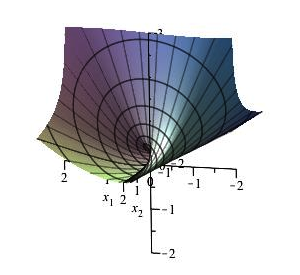}
\includegraphics[width=2in]{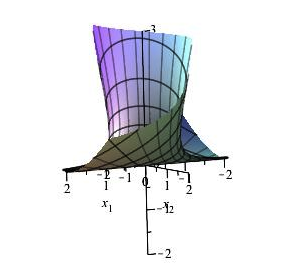} \includegraphics[width=2in]{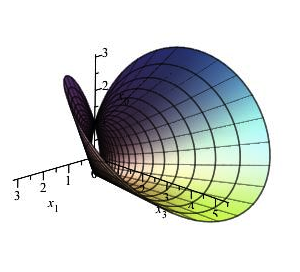}
\includegraphics[width=2in]{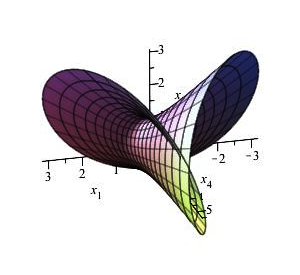}
\caption{Projections of surface from Example~\ref{time-like-example} to 3-D subspaces.}
\label{time-like-figure}
\end{figure}

As in the space-like examples, we can show that~$\Sigma$ is contained in the intersection of three quadric hypersurfaces
in $\mathbb{R}^5 \setminus \{0 \}$.
The coordinates of $\bar{f}(u,v)$ satisfy the quadratic equations
\begin{gather}
3 x_2^2 - x_4 (4 x_0 + 2) = 0,
\qquad
3 x_1^2 - x_3 (4 x_0 + 2) = 0,
\qquad
2 x_0^2 - x_0 - 3 x_1 x_2 - 1 = 0.
\label{time-like-implicit-eqns}
\end{gather}
Therefore,~$\Sigma$ is contained in the (real) algebraic variety $X \subset \mathbb{R}^5 \setminus \{0 \} $ def\/ined~by
equations~\eqref{time-like-implicit-eqns}.
$\Sigma$ is an interesting and somewhat complicated subset of~$X$: f\/irst
observe that the projection of~$X$ to the $(x_0, x_1, x_2)$ coordinate 3-plane consists of the hyperboloid of one sheet
def\/ined by the third equation in~\eqref{time-like-implicit-eqns}, minus all points of the form $(-\tfrac{1}{2}, x_1,
x_2)$ except for $(-\tfrac{1}{2}, 0,0)$.
The projection of~$X$ to this punctured hyperboloid is one-to-one, except over the point $(-\tfrac{1}{2}, 0,0)$, where
the inverse image consists of all points of the form $(-\tfrac{1}{2}, 0, 0, x_3, x_4)$.
Meanwhile, $\Sigma$ consists of that portion of~$X$ that projects to the portion of the hyperboloid with $x_0 >-\tfrac{1}{2}$,
together with the curve
\begin{gather*}
\big\{\big({-}\tfrac{1}{2}, 0, 0, x_3, x_4\big) \,|\, x_3 x_4 = \tfrac{9}{4}, \, x_3, x_4 > 0 \big\}.
\end{gather*}
Projections of~$\Sigma$ to the $(x_1, x_2, x_0)$, $(x_1, x_2, x_3)$, $(x_1, x_2, x_4)$, $(x_1, x_3, x_0)$, and $(x_1,
x_4, x_0)$ coordinate 3-planes are shown in Fig.~\ref{time-like-figure}.
\end{exmpl}

We collect the results of this section in the following theorem:

\begin{teo}
\label{homog-thm}
Let $\bar{f}:M \to \mathbb{R}^5 \setminus\{0\}$ be a~centroaffine immersion whose image $\Sigma = \bar{f}(M)$ is
a~homogeneous, nondegenerate, space-like or time-like centroaffine surface.
Then~$\Sigma$ is equivalent via the ${\rm GL}(5,\mathbb{R})$-action on $\mathbb{R}^5 \setminus \{0 \}$ to one of the
following:
\begin{itemize}\itemsep=0pt
\item the immersion of the hyperbolic plane $\mathbb{H}^2$ of Example~{\rm \ref{space-like-epsilon-plus-one}};
\item the immersion of the sphere~$S$ of Example~{\rm \ref{space-like-epsilon-minus-one}};
\item the immersion of the Lorentzian surface $S^{2,1}$ of Example~{\rm \ref{time-like-example}}.
\end{itemize}
\end{teo}

\subsection*{Acknowledgments}
This research was supported in part by NSF grants DMS-0908456 and DMS-1206272.

\pdfbookmark[1]{References}{ref}
\LastPageEnding

\end{document}